\documentstyle[10pt]{article}

\input amssym.def       
\input amssym.tex

\def\double{\Bbb}

\def\zzz{{\double Z}}

\def\ccc{{\double C}}
\def\dsum{\mathop{\sum}}

\newtheorem{lemma}{Lemma}[section]
\newtheorem{definition}{Definition}[section]
\newtheorem{proposition}{Proposition}[section]
\newtheorem{corollary}{Corollary}[section]
\newtheorem{remark}{Remark}[section]

\def\demo{\noindent\underline{Proof}:\par}
\def\edemo{\hfill$\square$\\}

\begin{document} 

\hsize 17truecm
\vsize 24truecm
\font\eightrm=cmr8
\baselineskip 14pt

\begin{titlepage}

\centerline{\Large CENTRE DE PHYSIQUE THEORIQUE}
\centerline{\large CNRS - Luminy, Case 907}
\centerline{\large 13288 Marseille Cedex 9}
\vskip 2truecm
\begin{center}
{\bf\Large  \sc $U_{q}(sl_{2})$ at fourth root of unity }
\end{center}
\bigskip
\begin{center}
{\bf D. Kastler, T. Krajewski, P. Seibt, K. Valavane}
\end{center}
\vskip 2truecm
\tableofcontents
\end{titlepage}

\section{Introduction}

This paper originated in a double motivation.

\par

In physics, we have strong reasons to believe that $SL_{q}(2,\ccc)$ for a primitive third root of unity is fundamentally related with the fermion structure \cite{AC1}. A precise implementation of this idea would open extremely interesting perspectives. In the noncommutative geometry version of the standard model of elementary particles \cite{CC1}\cite{CIKS}, replacement of the present phenomenological "spectral triple" by a fundamental Ansatz is expected to yield a scenario for a new "supersymmetric" standard model and a calculation of the fermion masses. 
 
\par

In mathematics, after Alain Connes axiomatic construction of spin-manifold \cite{AC1}, one is tempted by the project of an analogous theory of "generalized supermanifolds" (which we call "medusae").  This area, inevitable inasmuch as these objects exist (classical supermanifolds, gauge degrees of freedom of field theories) is mysterious and difficult, because the corresponding algebras are no longer subalgebras of $C^{*}$-algebras. In fact these are non semi-simple algebras, semi-simplicity arising only after quotienting by their nilradical. The projected survey of "medusae"lacks the two guiding features which have led Connes to the axiomatics of spin-manifold: the guidance through physics (hardly procured by the present supersymmetric standard model -in our opinion too ugly to be fundamental!) and (non semi-simplicity pending) the absence of Hilbert space techniques (presumably leading to abysses like "Hilbert spaces with indefinite metric", a situation already found but poor!
ly investigated in quantum field theory ).

\par

We believe that the (finite dimensional) $U_{q}(sl_{2})$ at roots of unity are models which hopefully yield features suggesting axioms for "medusae". A feature found in the two examples of  $U_{q}(sl_{2})$ for $q^{3}=1$ \cite{DK1} and the present $H_{1}^{i}$ is the striking (apparently new) fact that the trace of the adjoint representation (in the quantum group sense \cite{CASS}) has the nilradical in its kernel. This fact, in combination with an appropriate $*$-operation entails that semi-simplicity is synonymous with "positivity" (as is the case for the transversal degrees of freedom of quantum electrodynamics in the Lorentz gauge, eliminated by the requirement of a "strictly positive" Hilbert space, see e.g. \cite{DK2})

\par

We conjecture that semi-simplicity and positivity are synonymous for all roots of unity.

\par

Apart from this feature, our paper displays various aspects of the quantum groups $H^{i}_{N}$, including a complete description of their algebra automorphisms and Hopf $*$-structures.

\section{$U_{q}(sl_{2})$ at fourth root of unity}

$U_{q}(sl_{2})$ is the algebra defined by the symbols $K,K^{-1},E,F$ and
the relations 

\begin{equation}
\left\{ 
\begin{array}{c}
KE=q^{2}EK, \\ 
KF=q^{-2}FK, \\ 
\lbrack E,F]=\frac{K-K^{-1}}{q-q^{-1}}.
\end{array}
\right.  \label{1}
\end{equation}
It has a Hopf algebra structure defined by

\begin{equation}
Comultiplication\left\{ 
\begin{array}{c}
\Delta (\mathbf{1)=1\otimes 1}, \\ 
\Delta (K)\mathbf{=}K\otimes K,\\ 
\Delta (K^{-1})\mathbf{=}K^{-1}\otimes K^{-1}, \\ 
\Delta (E)\mathbf{=}E\otimes \mathbf{1}+K\otimes E, \\ 
\Delta (F)\mathbf{=}F\otimes K^{-1}+\mathbf{1}\otimes F,
\end{array}
\right.\label{2} 
\end{equation}

\begin{equation}
Counity\left\{ 
\begin{array}{c}
\varepsilon (\mathbf{1)=}\varepsilon (K)\mathbf{=}\varepsilon (K^{-1})=1, \\ 
\varepsilon (E)\mathbf{=}\varepsilon (F)=0,
\end{array}
\right.\label{2'}
\end{equation}

\begin{equation}
Antipode\left\{ 
\begin{array}{c}
S(\mathbf{1)=1\otimes 1}, \\ 
S(K)\mathbf{=}K^{-1}, \\ 
S(K^{-1})\mathbf{=}K, \\ 
S(E)\mathbf{=}-K^{-1}E, \\ 
S(F)\mathbf{=-}FK.
\end{array}
\right.  \label{3}
\end{equation}

When we consider the special case $q=i$, the relations $\ref{1}$ become

\begin{equation}
\left\{ 
\begin{array}{c}
KE=-EK, \\ 
KF=-FK, \\ 
\lbrack E,F]=\frac{K-K^{-1}}{2i}.
\end{array}
\right.  \label{4}
\end{equation}

The Casimir operator is defined as follows.

\begin{equation}
C=FE+\frac{K-K^{-1}}{4i}.  \label{5}
\end{equation}

\begin{lemma}
$E^{2},F^{2}$ and $K^{2}$ belong to the center of $U_{i}(sl_{2})$.
\end{lemma}

\demo
$K^{2}$ commutes with $K$, $E$ and $F$ since it anticommutes with $E,F$. $E^{2} $ commutes with $K$ since $E$ and $K$ anticommute. $E^{2}$ also
commutes with $F$, 

$$
\lbrack E^{2},F]=E[E,F]+[E,F]E=E\frac{K-K^{-1}}{2i}+\frac{K-K^{-1}}{2i}E=0. 
$$

Similarly, $F^{2}$ commutes with $K$ since $E$ and $K$ anticommute. $F^{2}$ commutes with $E$

$$
\lbrack E,F^{2}]=F[E,F]+[E,F]F=F\frac{K-K^{-1}}{2i}+\frac{K-K^{-1}}{2i}F=0.
$$
\edemo

\begin{definition}
We define $H_{N}^{i}=U_{i}(sl_{2})/I_{N}$, with $I_{N}$ the ideal of $U_{i}(sl_{2})$
generated by $E^{2},F^{2},K^{2N}-\mathbf{1.}$
\end{definition}

\begin{proposition}
1) $H_{N}^{i}$ is a Hopf algebra with the Hopf structure $\ref{2}$-$\ref{3}$ and
the Casimir operator $\ref{5}$.

2) $H_{N}^{i}$ has a PWB-base $\left\{ F^{p}K^{n}E^{q}\right\}
_{p,q=0,1;n=0,...,N}$. Thus it has dimension $8N.$
\end{proposition}

\demo
In fact, we have in $U_{i}(sl_{2})$

\begin{equation}
\left. 
\begin{array}{l}
\Delta (E^{2})=E^{2}\otimes \mathbf{1}+K^{2}\otimes E^{2}, \\ 
\Delta (F^{2})=F^{2}\otimes K^{-2}+\mathbf{1}\otimes F^{2}, \\ 
\Delta (K^{2N}-\mathbf{1})=K^{2N}\otimes (K^{2N}-\mathbf{1)+}(K^{2N}-\mathbf{%
1)}\otimes \mathbf{1},
\end{array}
\right.  \label{6}
\end{equation}

\begin{equation}
\left. 
\begin{array}{l}
\varepsilon (E^{2})=0, \\ 
\varepsilon (F^{2})=0, \\ 
\varepsilon (K^{2N}-\mathbf{1})=0,
\end{array}
\right. \left. 
\begin{array}{l}
S(E^{2})=-K^{-2}E^{2}, \\ 
S(F^{2})=-F^{2}K^{2}, \\ 
S(K^{2N}-\mathbf{1})=K^{-2N}(K^{2N}-\mathbf{1}),
\end{array}
\right.  \label{7}
\end{equation}

1) The relations are immediately checked using multiplicativity of $\Delta $ and $\varepsilon $, as well as antimultiplicativity of $S$. 
These relations imply that $I_{N}$ is a Hopf ideal. Indeed one has from $\ref
{6}$ the inclusion $\Delta (I_{N})\subset I_{N}\otimes
H_{N}^{i}+H_{N}^{i}\otimes I_{N}$, and $I_{N}$ contains the elements $\varepsilon (E^{2})$, $\varepsilon (F^{2})$, $\varepsilon (K^{2N}-\mathbf{1})$, $S(E^{2})$, $S(F^{2}),S(K^{2N}-\mathbf{1})$ owing to $\ref{7}$.

2) Follows from the multiplication table below.
\edemo

We now introduce a convenient alternative parametrization of the $N$%
-dimensional algebra $\mathbf{K}$ generated by $K$. Owing to $K^{2N}-\mathbf{1}$, $\mathbf{K}$ is the group algebra of the finite abelian group $\zzz/2N\zzz$. Using harmonic analysis on this group, we replace the basis $\{K^{n}\}_{n=0,..,2N-1}$ by the Fourier transformed basis $\{e_{n}\}_{n=0,..,2N-1}$, leading to a simpler description.

\begin{lemma}
Setting a complex number $u=e^{\frac{2i\pi }{2N}}$, the $*$-symmetric
elements : 
\begin{equation}
e_{k}=e_{k}^{*}=\frac{1}{2N}\dsum\limits_{j\in \zzz/2N\zzz}u^{kj}K^{j}\qquad k\in \zzz/2N\zzz,  \label{8}
\end{equation}

with reversion formulae

\begin{equation}
K^{j}=\dsum\limits_{k\in \zzz/2N\zzz}u^{-kj}e_{k}\qquad j\in \zzz/2N\zzz,  \label{9}
\end{equation}

yield a basis of $K$. It has the following properties, 
\begin{equation}
\left. 
\begin{array}{l}
1)\qquad \dsum\limits_{k\in \zzz/2N\zzz}e_{k}=\mathbf{1,} \\ 
2)\qquad e_{k}e_{m}=\delta _{km}e_{k}, \\ 
3)\qquad Ke_{k}=u^{-k}e_{k}, \\ 
4)\qquad K^{-1}e_{k}=u^{k}e_{k}, \\ 
5)\qquad Ee_{k}=e_{k+N}E, \\ 
6)\qquad e_{k}F=Fe_{k+N}, \\ 
7)\qquad C=EF+\frac{1}{4i}\dsum\limits_{k\in \zzz/2N\zzz%
}(u^{k}-u^{-k})e_{k}.
\end{array}
k\in \zzz/2N\zzz\right.  \label{10}
\end{equation}
\end{lemma}

\demo
The equivalence of $\ref{8}$ and $\ref{9}$ stems from the obvious fact that $%
\dsum\limits_{k\in \zzz/2N\zzz}u^{(k-m)j}=\delta _{km}$. Check of the
other claims: We have

\begin{eqnarray*}
2)\quad e_{k}e_{m} &=&\dsum\limits_{j\in \zzz/2N\zzz%
}u^{kj}K^{j}e_{m}=\left( \dsum\limits_{j\in \zzz/2N\zzz%
}u^{(k-m)j}\right) e_{m}=\delta _{km}e_{k} \\
3)\quad Ke_{k} &=&\frac{1}{2N}K\left( \dsum\limits_{j\in \zzz/2N\zzz%
}u^{kj}K^{j}\right) \\
&=&\frac{u^{-k}}{2N}\dsum\limits_{j\in \zzz/2N\zzz%
}u^{k(j+1)}K^{j+1}=u^{-k}e_{k} \\
5)\quad Ee_{k} &=&\frac{1}{2N}E\left( \dsum\limits_{j\in \zzz/2N\zzz%
}u^{kj}K^{j}\right) =\frac{1}{2N}\dsum\limits_{j\in \zzz/2N\zzz%
}u^{(N+k)j}K^{j}E=e_{k+N}E \\
6)\quad e_{k}F &=&\frac{1}{2N}\left( \dsum\limits_{j\in \zzz/2N\zzz%
}u^{kj}K^{j}\right) F=\frac{1}{2N}F\dsum\limits_{j\in \zzz/2N\zzz%
}u^{(N+k)j}K^{j}=Fe_{k+N}
\end{eqnarray*}

Let us now describe in detail the two cases $N=1$ and $N=2$. 

\subsection{Case $N=1$}

For $N=1$ we have $
K^{2}=\mathbf{1}$, hence $K=K^{-1}$.

\begin{definition}
The algebra $H_{1}^{i}$ is defined by the symbols $K,E,F$ and the relations
\begin{equation}
\left\{ 
\begin{array}{c}
KE=-EK, \\ 
KF=-FK, \\ 
\lbrack E,F]=0,
\end{array}
\right. \qquad \left\{ 
\begin{array}{c}
E^{2}=0, \\ 
F^{2}=0, \\ 
K^{2}=\mathbf{1}.
\end{array}
\right.  \label{11}
\end{equation}
\end{definition}

\begin{lemma}
1) With 
\begin{equation}
\left\{ 
\begin{array}{c}
e_{0}=\frac{\mathbf{1}+K}{2}, \\ 
e_{1}=\frac{\mathbf{1}-K}{2},
\end{array}
\right.  \label{12}
\end{equation}
$H_{1}^{i}$ is equivalently defined by the symbols $e_{0},e_{1},E,F,K$ and
the relations 
\begin{equation}
\left\{ 
\begin{array}{c}
e_{0}^{2}=e_{0}, \\ 
e_{1}^{2}=e_{1}, \\ 
e_{1}e_{0}=0, \\ 
e_{0}+e_{1}=\mathbf{1},
\end{array}
\right. \left\{ 
\begin{array}{c}
E^{2}=0, \\ 
F^{2}=0, \\ 
\lbrack E,F]=0,
\end{array}
\right. \left\{ 
\begin{array}{c}
e_{0}E=Ee_{1}, \\ 
e_{1}E=Ee_{0}, \\ 
e_{0}F=Fe_{1}, \\ 
e_{1}F=Fe_{0},
\end{array}
\right.  \label{13}
\end{equation}

yielding the basis 
\begin{equation}
\left\{ 
\begin{array}{c}
e_{0}, \\ 
e_{1}, \\ 
E_{0}=e_{0}E, \\ 
E_{1}=e_{1}E, \\ 
F_{0}=e_{0}F, \\ 
F_{1}=e_{1}F, \\ 
C_{0}=e_{0}EF, \\ 
C_{1}=e_{1}EF.
\end{array}
\right.  \label{14}
\end{equation}

2) $H_{1}^{i}$ is a $*$-algebra under the $*$-operation specified by 
$$
\left\{ 
\begin{array}{c}
E^{*}=F, \\ 
F^{*}=E, \\ 
K^{*}=K,
\end{array}
\right. \qquad \left\{ 
\begin{array}{c}
e_{0}^{*}=e_{0}, \\ 
e_{1}^{*}=e_{1}.
\end{array}
\right. 
$$

3) The Casimir element is
$$
C=C^{*}=EF=FE. 
$$

4) The basis $\ref{14}$ is acted upon as follows by multiplications by $E$, $F$
and $C$: 
\begin{equation}
\left\{ 
\begin{array}{c}
Ee_{0}=E_{1}, \\ 
Ee_{1}=E_{0}, \\ 
EE_{0}=0, \\ 
EE_{1}=0, \\ 
EF_{0}=C_{1}, \\ 
EF_{1}=C_{0}, \\ 
EC_{0}=0, \\ 
EC_{1}=0,
\end{array}
\right. \left\{ 
\begin{array}{c}
Fe_{0}=F_{1}, \\ 
Fe_{1}=F_{0}, \\ 
FE_{0}=C_{1}, \\ 
FE_{1}=C_{0}, \\ 
FF_{0}=0, \\ 
FF_{1}=0, \\ 
FC_{0}=0, \\ 
FC_{1}=0,
\end{array}
\right. \left\{ 
\begin{array}{c}
e_{0}E=E_{0}, \\ 
e_{1}E=E_{1}, \\ 
E_{0}E=0, \\ 
E_{1}E=0, \\ 
F_{0}E=C_{0}, \\ 
F_{1}E=C_{1}, \\ 
C_{0}E=0, \\ 
C_{1}E=0,
\end{array}
\right. \left\{ 
\begin{array}{c}
e_{0}F=F_{0}, \\ 
e_{1}F=F_{1}, \\ 
E_{0}F=C_{0}, \\ 
E_{1}F=C_{1}, \\ 
F_{0}F=0, \\ 
F_{1}F=0, \\ 
C_{0}F=0, \\ 
C_{1}F=0,
\end{array}
\right. \left\{ 
\begin{array}{c}
Ce_{0}=C_{0}, \\ 
Ce_{1}=C_{1}, \\ 
CE_{0}=0, \\ 
CE_{1}=0, \\ 
CF_{0}=0, \\ 
CF_{1}=0, \\ 
CC_{0}=0, \\ 
CC_{1}=0.
\end{array}
\right.  \label{15}
\end{equation}

5) The multiplication table of $H_{1}^{i}$ is as follows (we plugged the
product $XY$ at the intersection of line $X$ and column $Y$, the latter so as to
have $*$-symmetry w.r.t. the diagonal). 
\begin{center}
\begin{tabular}{|c|c|c|c|c|c|c|c|c|}
\hline
& $e_{0}$ & $e_{1}$ & $E_{0}$ & $E_{1}$ & $F_{0}$ & $F_{1}$ & $C_{0}$ & $C_{1}$ \\ 
\hline
$e_{0}$ & $e_{0}$ & $0$ & $E_{0}$ & $0$ & $F_{0}$ & $0$ & $C_{0}$ & $0$ \\ 
\hline
$e_{1}$ & $0$ & $e_{1}$ & $0$ & $E_{1}$ & $0$ & $F_{1}$ & $0$ & $C_{1}$ \\ 
\hline
$E_{0}$ & $0$ & $E_{0}$ & $0 $& $0$ & $0$ & $C_{0}$ & $0$ & $0$ \\ 
\hline
$E_{1}$ & $E_{1}$ & $0$ & $0$ & $0$ & $C_{1}$ & $0$ & $0$ & $0$ \\ 
\hline
$F_{0}$ & $0$ & $F_{0}$ & $0$ & $C_{0}$ & $0$ & $0$ & $0$ & $0$ \\ 
\hline
$F_{1}$ & $F_{1}$ & $0$ & $C_{1}$ & $0$ & $0$ & $0$ & $0$ & $0$ \\ 
\hline
$C_{0}$ & $C_{0}$ & $0$ & $0$ & $0$ & $0$ & $0$ & $0$ & $0$ \\ 
\hline
$C_{1}$ & $0$ & $C_{1}$ & $0$ & $0$ & $0$ & $0$ & $0$ & $0$\\
\hline
\end{tabular}
\end{center}
\begin{equation}
\label{16}
\end{equation}
\end{lemma}

\demo
Immediate from relations $\ref{11}$ and $\ref{13}$.
\edemo

\subsection{Structure  of the algebra $H_{1}^{i}$}

\begin{proposition}
1) Let $\mathbf{\ M}_{2}$ be algebra of the $2\times 2$ complex matrices,
equipped with its natural grading (diagonal entries are even and off-diagonal
entries odd) and $\Lambda =\Lambda _{1}\otimes \Lambda _{1}$ with $\Lambda
_{1}$ the Grassmann algebra over $\ccc$ (ordinary-not skew-tensor product
equipped with the tensor product $\zzz/2\zzz$-grading).

As an algebra $H_{1}^{i}$ is isomorphic to the even part $\left( 
\mathbf{M}_{2}\otimes \Lambda \right) ^{+}$ (again ordinary-not skew- tensor 
\rule{0.01in}{0.01in}product equipped with the tensor product $\zzz/2\Bbb{%
Z}$-grading). This isomorphism is specified as follows. With $\mathbf{M}%
_{2} $ spanned by its matrix units $\left\{ e_{lk}\right\} _{l,k=0,1}$; the
first $\Lambda _{1}$-factor by $\mathbf{1},e$; the second factor by $\mathbf{%
1},f$; and the tensor product $\Lambda _{1}\otimes \Lambda _{1}$ by $\mathbf{%
1}\otimes \mathbf{1},E=e\otimes \mathbf{1},F=\mathbf{1}\otimes f,EF=e\otimes
f$, one has,
$$
\left\{ 
\begin{array}{c}
e_{0}=e_{00}\otimes \mathbf{1}, \\ 
e_{1}=e_{11}\otimes \mathbf{1},
\end{array}
\begin{array}{c}
E_{0}=e_{01}\otimes E, \\ 
E_{1}=e_{10}\otimes E,
\end{array}
\begin{array}{c},
F_{0}=e_{00}\otimes F, \\ 
F_{1}=e_{11}\otimes F,
\end{array}
\begin{array}{c},
C_{0}=e_{00}\otimes EF, \\ 
C_{1}=e_{11}\otimes EF.
\end{array}
\right. 
$$

2) The subspace $N_{i}^{1}$ of $H_{1}^{i}$ spanned by $%
E_{0},E_{1},F_{0},F_{1},C_{0}$ and $C_{1}$ is the latter's nilradical,
giving rise to the quotient $H_{1}^{i}/N_{i}^{1}\cong\ccc\oplus\ccc$.

3) We definine as follows the scalar product $<.,.>$, 
$$
<a,b>=Tr\lambda (a^{*}b)\qquad a,b\in H_{1}^{i} 
$$

where $\lambda $ denotes the left-regular representation of $H_{1}^{i}$. $
<.,.>$ is positive semi-definite with null-space the nilradical $N_{i}^{1}$.
It is positive definite on the span of $e_{0},e_{1}$, where it coincides
with the usual trace of $\ccc\oplus\ccc$ ($%
<e_{l},e_{k}>=\delta _{lk}$, $l,k=0,1$).
\end{proposition}

\demo
1) One immediately checks that the elements \ref{14} fulfill the
multiplication rules $\ref{15}$.

2) These multiplication rules imply that $N_{i}^{1}$ is a subalgebra
fulfilling the inclusions $e_{0}N_{1}^{i},N_{1}^{i}e_{0},e_{1}N_{1}^{i},N_{1}^{i}e_{1}\subset N_{1}^{i}$, $N_{1}^{i}$ is thus an ideal of $H_{1}^{i}$ which moreover consists of nilpotent elements and yields a quotient
generated by $e_{0},e_{1}$ fulfilling $e_{0}^{2}=e_{0}$, $e_{1}^{2}=e_{1}$, $%
e_{0}e_{1}=0$, thus isomorphic to the semi-simple $\ccc\oplus\ccc$. Accordingly, $N_{i}^{1}$ is the nilradical of $H_{1}^{i}$.

3) Clear by inspection of $\ref{15}$, implying first that $Tr(\rho )$ vanishes
on $N_{i}^{1}$ and that $Tr\rho (e_{0})=Tr\rho (e_{1})=1$.
\edemo

\subsection{Case $N=2$}

We now describe the case corresponding to $N=2$, one has now $K^{4}=\mathbf{1%
}$ and $K-K^{-1}$ no longer vanishes.

\begin{definition}
The algebra $H_{2}^{i}$ is defined by the symbols $K$, $E$ and $F$ together with the relations
$$
\left\{ 
\begin{array}{c}
KE=-EK, \\ 
KF=-FK, \\ 
\lbrack E,F]=\frac{1}{2i}(K-K^{-1}),
\end{array}
\right. \qquad \left\{ 
\begin{array}{c}
E^{2}=0, \\ 
F^{2}=0, \\ 
K^{4}=\mathbf{1}.
\end{array}
\right. 
$$
\end{definition}

\begin{lemma}
1) Let $e_{m}=\frac{1}{4}\dsum\limits_{k\in \zzz/4\zzz}i^{mk}K^{k}$,
i.e 
$$
\left\{ 
\begin{array}{c}
e_{0}=\frac{1}{4}\left( \mathbf{1}+K+K^{2}+K^{3}\right), \\ 
e_{1}=\frac{1}{4}\left( \mathbf{1}+iK-K^{2}-iK^{3}\right), \\ 
e_{2}=\frac{1}{4}\left( \mathbf{1}-K+K^{2}-K^{3}\right), \\ 
e_{3}=\frac{1}{4}\left( \mathbf{1}-K-K^{2}+iK^{3}\right).
\end{array}
\right. 
$$

It implies the relations 
$$
\left\{ 
\begin{array}{c}
K=e_{0}-ie_{1}-e_{2}+ie_{3}, \\ 
K^{-1}=e_{0}+ie_{1}-e_{2}-ie_{3},\\
K-K^{-1}=2i(e_{3}-e_{1}), 
\end{array}
\right.
$$

yielding the basis 
\begin{equation}
\left\{ 
\begin{array}{c}
e_{0}, \\ 
e_{2}, \\ 
E_{0}=e_{0}E, \\ 
E_{2}=e_{2}E, \\ 
F_{0}=e_{0}F, \\ 
F_{2}=e_{2}F, \\ 
P_{0}=C_{0}=e_{0}EF, \\ 
P_{2}=C_{2}=e_{2}EF,
\end{array}
\right. \qquad \left\{ 
\begin{array}{c}
e_{1}, \\ 
e_{3}, \\ 
E_{1}=e_{1}E, \\ 
E_{3}=e_{3}E, \\ 
F_{1}=e_{1}F, \\ 
F_{3}=e_{3}F, \\ 
P_{1}=C_{1}-\frac{1}{2}e_{1}=e_{1}EF, \\ 
P_{3}=C_{3}+\frac{1}{2}e_{3}=e_{3}EF.
\end{array}
\right.  \label{17}
\end{equation}

Accordingly, $H_{2}^{i}$ is equivalently defined by the symbols $E$, $F$ and $e_{j}$ $j=0,1,2,3$ and the relations  
$$
\left\{ 
\begin{array}{c}
e_{0}+e_{1}+e_{2}+e_{3}=\mathbf{1}, \\ 
e_{l}e_{k}=\delta _{lk},
\end{array}
,l=0,1,2,3\right. \left\{ 
\begin{array}{c}
E^{2}=0, \\ 
F^{2}=0, \\ 
\lbrack E,F]=e_{3}-e_{1},
\end{array}
\right. \left\{ 
\begin{array}{c}
e_{j}E=Ee_{j+2}, \\ 
Fe_{j}=e_{j+2}F,
\end{array}
j\in \zzz/4\zzz.\right. 
$$

2) $H_{2}^{i}$ is a $*$-algebra under the $*$-operation specified by 
$$
\left\{ 
\begin{array}{c}
E^{*}=F, \\ 
F^{*}=E,
\end{array}
\right. \qquad \left\{ 
\begin{array}{c}
K^{*}=K, \\ 
e_{j}^{*}=e_{j},
\end{array}
\right. j=0,1,2,3. 
$$

3) The Casimir element is 
\begin{equation}
C=C^{*}=FE-\frac{1}{2}(e_{1}-e_{3})=EF+\frac{1}{2}(e_{1}-e_{3}).  \label{18}
\end{equation}

4) Let 
$$
\left\{ 
\begin{array}{c}
\pi _{0}=e_{0}+e_{2}, \\ 
\pi _{1}=e_{1}+e_{3},
\end{array}
,\right. \left\{ 
\begin{array}{c}
H_{2}^{i(0)}=\pi _{0}H_{2}^{i}, \\ 
H_{2}^{i(1)}=\pi _{1}H_{2}^{i}.
\end{array}
\right. 
$$

Then $\pi _{0}$ and $\pi _{1}$ are supplementary central idempotents
yielding supplementary ideals $H_{2}^{i(0)}$ and $H_{2}^{i(1)}$ of $
H_{2}^{i} $ spanned respectively by the left and right basis $\ref{17}$.

5) Multiplication by $E$, $F$ and $C$ act as follows on the basis $\ref{17}$, 
$$
\left\{ 
\begin{array}{c}
Ee_{0}=E_{2}, \\ 
Ee_{2}=E_{0}, \\ 
EE_{0}=0, \\ 
EE_{2}=0, \\ 
EF_{0}=P_{2}, \\ 
EF_{2}=P_{0}, \\ 
EP_{0}=0, \\ 
EP_{1}=0,
\end{array}
\right. \left\{ 
\begin{array}{c}
Fe_{0}=F_{2}, \\ 
Fe_{2}=F_{0}, \\ 
FE_{0}=P_{2}, \\ 
FE_{2}=P_{0}, \\ 
FF_{0}=0, \\ 
FF_{2}=0, \\ 
FP_{0}=0, \\ 
FP_{1}=0.
\end{array}
\right. \left\{ 
\begin{array}{c}
e_{0}E=E_{0}, \\ 
e_{2}E=E_{2}, \\ 
E_{0}E=0, \\ 
E_{2}E=0, \\ 
F_{0}E=P_{0}, \\ 
F_{2}E=P_{2}, \\ 
P_{0}E=0, \\ 
P_{2}E=0,
\end{array}
\right. \left\{ 
\begin{array}{c}
e_{0}F=F_{0}, \\ 
e_{2}F=F_{2}, \\ 
E_{0}F=P_{0}, \\ 
E_{2}F=P_{2}, \\ 
F_{0}F=0, \\ 
F_{2}F=0, \\ 
P_{0}F=0, \\ 
P_{1}F=0,
\end{array}
\right. \left\{ 
\begin{array}{c}
Pe_{0}=P_{0}, \\ 
Pe_{2}=P_{2}, \\ 
PE_{0}=0, \\ 
PE_{2}=0, \\ 
PF_{0}=0, \\ 
PF_{2}=0, \\ 
PP_{0}=0, \\ 
PP_{2}=0,
\end{array}
\right. 
$$
and 
$$
\left\{ 
\begin{array}{c}
Ee_{1}=E_{3}, \\ 
Ee_{3}=E_{1}, \\ 
EE_{1}=0, \\ 
EE_{3}=0, \\ 
EF_{1}=P_{3}, \\ 
EF_{3}=P_{1}, \\ 
EP_{1}=0, \\ 
EP_{3}=0,
\end{array}
\right. \left\{ 
\begin{array}{c}
Fe_{1}=F_{3}, \\ 
Fe_{3}=F_{1}, \\ 
FE_{1}=P_{3}, \\ 
FE_{3}=P_{1}, \\ 
FF_{1}=0, \\ 
FF_{3}=0, \\ 
FP_{1}=0, \\ 
FP_{3}=0,
\end{array}
\right. \left\{ 
\begin{array}{c}
e_{1}E=E_{1}, \\ 
e_{3}E=E_{3}, \\ 
E_{1}E=0, \\ 
E_{3}E=0, \\ 
F_{1}E=P_{1}+e_{1}, \\ 
F_{3}E=P_{3}-e_{3}, \\ 
P_{1}E=-E_{1}, \\ 
P_{3}E=E_{3},
\end{array}
\right. \left\{ 
\begin{array}{c}
e_{1}F=F_{1}, \\ 
e_{3}F=F_{3}, \\ 
E_{1}F=P_{1}, \\ 
E_{3}F=P_{3}, \\ 
F_{1}F=0, \\ 
F_{3}F=0, \\ 
P_{1}F=0, \\ 
P_{3}F=0,
\end{array}
\right. \left\{ 
\begin{array}{c}
Ce_{1}=C_{1}, \\ 
Ce_{3}=C_{3}, \\ 
CE_{1}=-\frac{1}{2}E_{1}, \\ 
CE_{3}=\frac{1}{2}E_{3}, \\ 
CF_{1}=\frac{1}{2}F_{1}, \\ 
CF_{3}=-\frac{1}{2}F_{3}, \\ 
CP_{1}=-\frac{1}{2}P_{1}, \\ 
CP_{3}=\frac{1}{2}P_{3}.
\end{array}
\right. 
$$
6) One has the following multiplication tables (the product $XY$ is at the
intersection of line $X$ and column $Y$, ordering the latter so as to have * symmetry w.r.t the diagonal) for $H_{2}^{i(0)}$ in the basis $\left\{
e_{0},e_{2},E_{0},E_{2},F_{0},F_{2},P_{0},P_{2}\right\}$
\begin{center}
\begin{tabular}{|c|c|c|c|c|c|c|c|c|}
\hline
& $e_{0}$ & $e_{2}$ & $E_{0}$ & $E_{2}$ & $F_{0}$ & $F_{2}$ & $P_{0}$ & $P_{2}$ \\ 
\hline
$e_{0}$ & $e_{0}$ & $0$ & $E_{0}$ & $0$ & $F_{0}$ & $0$ & $P_{0}$ & $0$ \\ 
\hline
$e_{2}$ & $0$ & $e_{2}$ & $0$ & $E_{2}$ & $0$ & $F_{2}$ & $0$ & $P_{2}$ \\ 
\hline
$E_{0}$ & $0$ & $E_{0}$ & $0$ & $0$ & $0$ & $P_{0}$ & $0$ & $0$ \\ 
\hline
$E_{2}$ & $E_{2}$ & $0$ & $0$ & $0$ & $P_{2}$ & $0$ & $0$ & $0$ \\ 
\hline
$F_{0}$ & $0$ & $F_{0}$ & $0$ & $P_{0}$ & $0$ & $0$ & $0$ & $0$ \\ 
\hline
$F_{2}$ & $F_{2}$ & $0$ & $P_{2}$ & $0$ & $0$ & $0$ & $0$ & $0$ \\ 
\hline
$P_{0}$ & $P_{0}$ & $0$ & $0$ & $0$ & $0$ & $0$ & $0$ & $0$ \\ 
\hline
$P_{2}$ & $0$ & $P_{2}$ & $0$ & $0$ & $0$ & $0$ & $0$ & $0$\\
\hline
\end{tabular}
\end{center}
\begin{equation}
\label{19}
\end{equation}

and for $H_{2}^{i(1)}$ in the basis $\left\{
e_{1},e_{3},E_{1},E_{3},F_{1},F_{3},P_{1},P_{3}\right\} $
\begin{center}
\begin{tabular}{|c|c|c|c|c|c|c|c|c|}
\hline
& $e_{1}$ & $e_{3}$ & $E_{1}$ & $E_{3}$ & $F_{1}$ & $F_{3}$ & $P_{1}$ & $P_{3}$ \\ 
\hline
$e_{1}$ & $e_{1}$ & $0$ & $E_{1}$ & $0$ & $F_{1}$ & $0$ & $P_{1}$ & $0$ \\ 
\hline
$e_{3}$ & $0$ & $e_{3}$ & $0$ & $E_{3}$ & $0$ & $F_{3}$ & $0$ & $P_{3}$ \\ 
\hline
$E_{1}$ & $0$ & $E_{1}$ & $0$ & $0 $& $0 $& $P_{1}$ & $0$ & $0$ \\ 
\hline
$E_{3}$ & $E_{3}$ & $0$ & $0$ & $0$ & $P_{3}$ & $0$ & $0$ & $0$ \\ 
\hline
$F_{1}$ & $0$ & $F_{1}$ & $0$ & $P_{1}+e_{1}$ & $0$ & $0$ & $0$ & $F_{1}$ \\ 
\hline
$F_{3}$ & $F_{3}$ & $0$ & $P_{3}-e_{3}$ & $0$ & $0$ & $0$ & $-F_{3}$ & $0$ \\ 
\hline
$P_{1}$ & $P_{1}$ & $0$ & $-E_{1}$ & $0$ & $0$ & $0$ & $-P_{1}$ & $0$ \\ 
\hline
$P_{3}$ & $0$ & $P_{3}$ &$ 0 $& $E_{3}$ & $0$ & $0$ & $0$ & $P_{3}$\\
\hline
\end{tabular}
\end{center}
\begin{equation}
\label{20}
\end{equation}

This yields the following action of the Casimir operator 
\begin{center}
\begin{tabular}{|c|c|c|c|c|c|c|c|c|}
\hline
$x$& $e_{1}$ & $e_{3}$ & $E_{1}$ & $E_{3}$ & $F_{1}$ & $F_{3}$ & $P_{1}$ & $P_{3}$ \\ 
\hline
$C(x)$ & $C_{1}$ & $C_{3}$ & $-\frac{1}{2}E_{1}$ & $\frac{1}{2}E_{3}$ & $\frac{1}{2}F_{1}$ & $-\frac{1}{2}E_{3}$ & $-\frac{1}{2}P_{1}$ & $\frac{1}{2}P_{3}$\\
\hline
\end{tabular}
\end{center}
and shows that the restriction of $C$ to $H_{2}^{i(1)}$ has the eigenspaces $
\ccc E_{1}+\ccc F_{3}+\ccc P_{1}$ to the eigenvalue $-\frac{1}{2}$ and 
$\ccc E_{3}+\ccc F_{1}+\ccc P_{3}$ to the eigenvalue $\frac{1}{2}$.
\end{lemma}

\demo
The products $\ref{19}$ and $\ref{20}$ not involving the $C_{j}$
are immediate from $\ref{11}$ and  $\ref{13}$. Check of the products $\ref{20}$ involving $P_{1},P_{3}$ is made using

$$
\left\{ 
\begin{array}{l}
EFE=-(e_{1}-e_{3})E, \\ 
FEF=(e_{1}-e_{3})F,
\end{array}
\right. 
$$

so that, 

$$
\left\{ 
\begin{array}{c}
E_{1}F_{3}=e_{1}EF=P_{1}, \\ 
F_{3}E_{1}=e_{3}FE=e_{3}(EF+e_{1}-e_{3})=P_{3}-e_{3},
\end{array}
\right. 
$$

$$
\left\{ 
\begin{array}{l}
F_{3}P_{1}=e_{3}FEF=e_{3}(e_{1}-e_{3})F=-F_{3},\\
P_{1}E_{1}=e_{1}EFE=-e_{1}(e_{1}-e_{3})E=-E_{1},\\
P_{1}P_{1}=e_{1}EFEF=-e_{1}(e_{1}-e_{3})EF=-e_{1}EF,
\end{array}
\right. 
$$

$$
\left\{ 
\begin{array}{l}
E_{3}F_{1}=e_{3}EF=P_{3},\\
F_{1}E_{3}=e_{1}FE=e_{1}(EF+e_{1}-e_{3})=P_{1}+e_{1},\\
F_{1}P_{3}=e_{1}FEF=e_{1}(e_{1}-e_{3})F=e_{1}F=F_{1},
\end{array}
\right. 
$$

$$
\left\{ 
\begin{array}{l}
P_{3}E_{1}=e_{3}EFE=-e_{3}(e_{1}-e_{3})E=e_{3}E=E_{3},\\
P_{3}E_{3}=e_{3}EFE=-e_{3}(e_{1}-e_{3})E=e_{3}E=E_{3},\\
P_{3}P_{3}=e_{3}EFEF=-e_{3}(e_{1}-e_{3})EF=e_{3}EF=P_{3}.
\end{array}
\right. 
$$

\edemo

\subsection{Structure of the algebra $H_{2}^{i}$}

\begin{proposition}
1) The algebra $H_{2}^{i}$ splits into the direct sum of the ideals  $\footnote{algebraic (not Hopf) ideals.}H_{2}^{i(0)}$ and $H_{2}^{i(1)}$.

2) The algebra $H_{2}^{i(0)}$ is isomorphic to $H_{1}^{i}$ with the
isomorphism given by 
\begin{equation}
\left\{ 
\begin{array}{c}
e_{0}\rightarrow e_{0}, \\ 
e_{2}\rightarrow e_{1},
\end{array}
\right. \left\{ 
\begin{array}{c}
E_{0}\rightarrow E_{0}, \\ 
E_{2}\rightarrow E_{1},
\end{array}
\right. \left\{ 
\begin{array}{c}
F_{0}\rightarrow F_{0}, \\ 
F_{2}\rightarrow F_{1},
\end{array}
\right. \left\{ 
\begin{array}{c}
C_{0}\rightarrow C_{0}, \\ 
C_{2}\rightarrow C_{1}.
\end{array}
\right.  \label{21}
\end{equation}

Consequently, $H_{2}^{i(0)}$ is an algebra isomorphic to the even part $%
\left( \mathbf{M}_{2}\otimes \Lambda_{1}\otimes\Lambda_{1}\right) ^{+}$, the isomorphism being
specified as follows : with $\mathbf{M}_{2}$ spanned by its matrix units $%
\left\{ e_{lk}\right\} _{l,k=0,1}$; the first $\Lambda _{1}$-factor by $%
\mathbf{1}$, $e$, the second $\Lambda _{1}$-factor by $\mathbf{1}$, $f$; and
the tensor product $\Lambda _{1}\otimes \Lambda _{1}$ by $\mathbf{1}\otimes 
\mathbf{1},E=e\otimes \mathbf{1},F=\mathbf{1}\otimes f,EF=e\otimes f$, one
has 
$$
\left\{ 
\begin{array}{c}
e_{1}=e_{00}\otimes \mathbf{1}, \\ 
e_{2}=e_{11}\otimes \mathbf{1},
\end{array}
\begin{array}{c}
E_{0}=e_{01}\otimes E, \\ 
E_{2}=e_{10}\otimes E,
\end{array}
\begin{array}{c}
F_{0}=e_{01}\otimes F, \\ 
F_{2}=e_{10}\otimes F,
\end{array}
\begin{array}{c}
P_{0}=e_{00}\otimes EF, \\ 
P_{2}=e_{11}\otimes EF.
\end{array}
\right. 
$$

3) The nilradical $N_{1}^{i(0)}$ of $H_{2}^{i(0)}$ is the eigenspace of the
Casimir element $C$ to the eigenvalue $0$, it is spanned by the elements $%
E_{0},E_{2},F_{0},F_{2},P_{0},P_{2}$. The quotient algebra $H_{2}^{i(0)}/N_{1}^{i(0)}$ is isomorphic to $\ccc\oplus \ccc$.

4)$H_{2}^{i(1)}$ is an algebra isomorphic to the semi-simple algebra $M(2,\ccc)\oplus M(2,\ccc)$ with the isomorphism given by
$$
K=i\pmatrix{-1&0&0&0\cr0&1&0&0\cr0&0&-1&0\cr0&0&0&1},\;
E=\pmatrix{0&-1&0&0\cr0&0&0&0\cr0&0&0&0\cr0&0&1&0},\;
F=\pmatrix{0&0&0&0\cr1&0&0&0\cr0&0&0&1\cr0&0&0&0}.
$$

5)We define the scalar product $<.,.>$, 
$$
<a,b>=Tr\lambda (a^{*}b), 
$$

where $\lambda$ denotes the left regular representation of $H_{2}^{i}$. One
has that

\qquad a) the even and odd parts $H_{2}^{i(0)}$ and $H_{2}^{i(1)}$ are
mutually $<.,.>$-orthogonal,

\qquad b) its restriction to $H_{2}^{i(0)}$, $<.,.>$ is positive
semi-definite with null space $N_{1}^{i(0)}$, and is positive definite on
the span of $e_{0},e_{2}$, where it coincides with the usual trace of $\ccc\oplus \ccc$,

\qquad c) its restriction to $H_{2}^{i(1)}$, $<.,.>$ behaves as follows with 
\begin{eqnarray*}
H_{2}^{i(1,+)} &=&\;span\;of\;\left\{ \frac{1}{2}(e_{1}-P_{1}),\frac{1}{2}%
(e_{3}+P_{3}),F_{1}\right\}, \\
H_{2}^{i(1,0)} &=&\;span\;of\;\left\{ \frac{1}{2}(e_{1}+P_{1}),\frac{1}{2}%
(e_{3}-P_{3}),E_{1},E_{3}\right\}, \\
H_{2}^{i(1,-)} &=&\;span\;of\;\left\{ F_{3}\right\},
\end{eqnarray*}

$<.,.>$ is positive definite on $H_{2}^{i(1,+)}$, with $\frac{1}{2}%
(e_{1}-P_{1}),\frac{1}{2}(e_{3}+P_{3}),F_{1}$ orthonormal,

$<.,.>$ has $H_{2}^{i(1,0)}$ as its null space (space of vectors orthogonal
to all vectors),

$<.,.>$ is negative definite on $H_{2}^{i(1,-)}$, with $<F_{3},F_{3}>=-1$.
\end{proposition}

\demo
1) Recalls a former result, cf lemma $\ref{10}.2$.

2) The changes $\ref{21}$ turn $\ref{19}$ into $\ref{16}$.

3) The eigenspaces of $C$ acting on $H_{2}^{i(0)}$ are immediately found
computing $\ref{18}$ via $\ref{19}$. The (a priori known) fact that these subspaces
are ideals and the corresponding quotient are patent from $\ref{19}$ and $\ref{20}$.

4)Follows from the fact that these matrices satisfy the same relations and generate an 8-dimensional algebra. 

5) Inspection of $\ref{20}$ yields the following table of values of $Tr\lambda
(.)$ on $H_{2}^{i(1)}$

\begin{center}
\begin{tabular}{|c|c|c|c|c|c|c|c|c|}
\hline
$u$ & $e_{1}$ & $e_{3}$ & $E_{1}$ & $E_{3}$ & $F_{1}$ & $F_{3}$ & $P_{1}$ &
$P_{3}$ \\
\hline 
$Tr\lambda (u)$ & $1$ & $1$ & $0$ & $0$ & $1$ & $0$ & $0$ & $-1$\\
\hline
\end{tabular}
\end{center}

Thus we have the following table of values of the scalar product $<.,.>$ on $H_{2}^{i(1)}$ (we plotted $<u,v>$ at the intersection of line $u$ and $v$
\begin{center}
\begin{tabular}{|c|c|c|c|c|c|c|c|c|}
\hline
& $e_{1}$ & $e_{3}$ & $E_{1}$ & $E_{3}$ & $F_{1}$ & $F_{3}$ & $P_{1}$ & $P_{3}$ \\ 
\hline
$e_{1}$ & $1$ & $0$ & $0$ & $0$ & $0$ & $0$ & $-1$ & $0$ \\ 
\hline
$e_{3}$ & 0 & 1 & 0 & 0 & 0 & 0 & 0 & 1 \\ 
\hline
$E_{1}$ & 0 & 0 & 0 & 0 & 0 & 0 & 0 & 0 \\ 
\hline
$E_{3}$ & 0 & 0 & 0 & 0 & 0 & 0 & 0 & 0 \\ 
\hline
$F_{1}$ & 0 & 0 & 0 & 0 & 1 & 0 & 0 & 0 \\ 
\hline
$F_{3}$ & 0 & 0 & 0 & 0 & 0 & -1 & 0 & 0 \\ 
\hline
$P_{1}$ & -1 & 0 & 0 & 0 & 0 & 0 & 1 & 0 \\ 
\hline
$P_{3}$ & 0 & 1 & 0 & 0 & 0 & 0 & 0 & 1\\
\hline
\end{tabular}
\end{center}
The statement then immediately follows.
\edemo

\section{Hopf structure of $H_{1}^{i}$ and $H_{2}^{i}$}

Note that $\ref{2}$-$\ref{3}$ give the Hopf structure of $H_{1}^{i}$ and $H_{2}^{i}$. Recalling basic definitions, the algebra $H_{1}^{i}$ is
defined by symbols $K,E,F$ and the relations

$$
\left\{ 
\begin{array}{c}
KE=-EK, \\ 
KF=-FK, \\ 
\lbrack E,F]=0,
\end{array}
\right. \left\{ 
\begin{array}{c}
E^{2}=0, \\ 
F^{2}=0, \\ 
K^{2}=\mathbf{1},
\end{array}
\right. \left\{ 
\begin{array}{c}
e_{l}e_{k}=\delta _{lk},\quad l,k=0,1 \\ 
e_{0}+e_{1}=\mathbf{1},
\end{array}
\right. \left\{ 
\begin{array}{c}
e_{j}E=Ee_{j+1},\quad j\in \zzz/2\zzz \\ 
e_{j}F=Fe_{j+1}.\quad j\in \zzz/2\zzz
\end{array}
\right. 
$$
It is spanned by the basis

$$
e_{0}=\frac{\mathbf{1}+K}{2},e_{1}=\frac{\mathbf{1}-K}{2}%
,E_{0}=e_{0}E,E_{1}=e_{1}E,F_{0}=e_{0}F,C_{0}=e_{0}EF,C_{1}=e_{1}EF. 
$$

\subsection{Hopf structure of $H_{1}^{i}$}

\begin{proposition}
The coproduct is given by
\begin{equation}
\left\{ 
\begin{array}{c}
\Delta (e_{0})=e_{0}\otimes e_{0}+e_{1}\otimes e_{1}, \\ 
\Delta (e_{1})=e_{0}\otimes e_{1}+e_{1}\otimes e_{0}, \\ 
\Delta (E_{0})=E_{0}\otimes e_{0}+E_{1}\otimes e_{1}+e_{0}\otimes
E_{0}-e_{1}\otimes E_{1}, \\ 
\Delta (E_{1})=E_{0}\otimes e_{1}+E_{1}\otimes e_{0}+e_{0}\otimes
E_{1}-e_{1}\otimes E_{0}, \\ 
\Delta (F_{0})=F_{0}\otimes e_{0}-F_{1}\otimes e_{1}+e_{0}\otimes
F_{0}+e_{1}\otimes F_{1}, \\ 
\Delta (F_{1})=-F_{0}\otimes e_{1}+F_{1}\otimes e_{0}+e_{0}\otimes
F_{1}+e_{1}\otimes F_{0}, \\ 
\Delta (C_{0})=C_{0}\otimes e_{0}-C_{1}\otimes e_{1}+e_{0}\otimes
C_{0}-e_{1}\otimes C_{1}+E_{0}\otimes F_{0}+E_{1}\otimes F_{1}-F_{0}\otimes
E_{0}-F_{1}\otimes E_{1}, \\ 
\Delta (C_{1})=-C_{0}\otimes e_{1}+C_{1}\otimes e_{0}+e_{0}\otimes
C_{1}-e_{1}\otimes C_{0}+F_{0}\otimes E_{1}+F_{1}\otimes E_{0}+E_{0}\otimes
F_{1}+E_{1}\otimes F_{0}.
\end{array}
\right.  \label{22}
\end{equation}
The antipde and the counit are
\begin{equation}
\left\{ 
\begin{array}{c}
S(e_{0})=e_{0}, \\ 
S(e_{1})=e_{1}, \\ 
S(E_{0})=E_{1}, \\ 
S(E_{1})=-E_{0},
\end{array}
\right. \left\{ 
\begin{array}{c}
S(F_{0})=-F_{1}, \\ 
S(F_{1})=F_{0}, \\ 
S(C_{0})=C_{0}, \\ 
S(C_{1})=C_{1},
\end{array}
\right.  \label{23}
\end{equation}
\begin{equation}
\left\{ 
\begin{array}{c}
\varepsilon (e_{0})=1, \\ 
\varepsilon (e_{1})=\varepsilon (E_{0})=\varepsilon (E_{1})=\varepsilon
(F_{0})=\varepsilon (F_{1})=\varepsilon (C_{0})=\varepsilon (C_{1})=0.
\end{array}
\right.  \label{24}
\end{equation}
\end{proposition}

\demo
Check of $\ref{22}$ : Using $\ref{12}$, whence $K=e_{0}-e_{1}$, we have

\begin{eqnarray*}
\Delta (e_{0}) &=&\frac{1}{2}\Delta (\mathbf{1}+K)=\frac{1}{2}\left[ \mathbf{%
1\otimes 1}+K\mathbf{\otimes }K\right] =\frac{1}{2}\left[ (e_{0}+e_{1})%
\mathbf{\otimes }(e_{0}+e_{1})+(e_{0}-e_{1})\mathbf{\otimes }%
(e_{0}-e_{1})\right] \\
&=&e_{0}\otimes e_{0}+e_{1}\otimes e_{1},
\end{eqnarray*}

\begin{eqnarray*}
\Delta (e_{1}) &=&\frac{1}{2}\Delta (\mathbf{1}-K)=\frac{1}{2}\left[ \mathbf{%
1\otimes 1}-K\mathbf{\otimes }K\right] =\frac{1}{2}\left[ (e_{0}+e_{1})%
\mathbf{\otimes }(e_{0}+e_{1})-(e_{0}-e_{1})\mathbf{\otimes }%
(e_{0}-e_{1})\right] \\
&=&e_{0}\otimes e_{1}+e_{1}\otimes e_{0}
\end{eqnarray*}

\begin{eqnarray*}
\Delta (E_{0}) &=&\Delta (e_{0}E)=\Delta (e_{0})\Delta (E)=\left(
e_{0}\otimes e_{0}+e_{1}\otimes e_{1}\right) \left( E\otimes
(e_{0}+e_{1})+\left( e_{0}-e_{1}\right) \otimes E\right) \\
&=&E_{0}\otimes e_{0}+E_{1}\otimes e_{1}+e_{0}\otimes E_{0}-e_{1}\otimes
E_{1},
\end{eqnarray*}

\begin{eqnarray*}
\Delta (E_{1}) &=&\Delta (e_{1}E)=\Delta (e_{1})\Delta (E)=\left(
e_{0}\otimes e_{1}+e_{1}\otimes e_{0}\right) \left( E\otimes
(e_{0}+e_{1})+\left( e_{0}-e_{1}\right) \otimes E\right) \\
&=&E_{0}\otimes e_{1}+E_{1}\otimes e_{0}+e_{0}\otimes E_{1}-e_{1}\otimes
E_{0}.
\end{eqnarray*}

Further, taking account of the fact that $K^{-1}=K=e_{0}-e_{1}$

\begin{eqnarray*}
\Delta (F_{0}) &=&\Delta (e_{0}F)=\Delta (e_{0})\Delta (F)=\left(
e_{0}\otimes e_{0}+e_{1}\otimes e_{1}\right) \left( F\otimes \left(
e_{0}-e_{1}\right) +(e_{0}+e_{1})\otimes F\right) \\
&=&F_{0}\otimes e_{0}-F_{1}\otimes e_{1}+e_{0}\otimes F_{0}+e_{1}\otimes
F_{1},
\end{eqnarray*}

\begin{eqnarray*}
\Delta (F_{1}) &=&\Delta (e_{1}F)=\Delta (e_{1})\Delta (F)=\left(
e_{0}\otimes e_{1}+e_{1}\otimes e_{0}\right) \left( F\otimes \left(
e_{0}-e_{1}\right) +(e_{0}+e_{1})\otimes F\right) \\
&=&-F_{0}\otimes e_{1}+F_{1}\otimes e_{0}+e_{0}\otimes F_{1}+e_{1}\otimes
F_{0},
\end{eqnarray*}

\begin{eqnarray*}
\Delta (C_{0}) &=&\Delta (E_{0})\Delta (F) \\
&=&\left( E_{0}\otimes e_{0}+E_{1}\otimes e_{1}+e_{0}\otimes
E_{0}-e_{1}\otimes E_{1}\right) \left( F\otimes \left( e_{0}-e_{1}\right)
+(e_{0}+e_{1})\otimes F\right) \\
&=&C_{0}\otimes e_{0}-C_{1}\otimes e_{1}+e_{0}\otimes C_{0}-e_{1}\otimes
C_{1}+E_{0}\otimes F_{0}+E_{1}\otimes F_{1}-F_{0}\otimes E_{0}-F_{1}\otimes
E_{1},
\end{eqnarray*}

\begin{eqnarray*}
\Delta (C_{1}) &=&\Delta (E_{1})\Delta (F) \\
&=&\left( E_{0}\otimes e_{1}+E_{1}\otimes e_{0}+e_{0}\otimes
E_{1}-e_{1}\otimes E_{0}\right) \left( F\otimes \left( e_{0}-e_{1}\right)
+(e_{0}+e_{1})\otimes F\right) \\
&=&-C_{0}\otimes e_{1}+C_{1}\otimes e_{0}+e_{0}\otimes C_{1}-e_{1}\otimes
C_{0}+F_{0}\otimes E_{1}+F_{1}\otimes E_{0}+E_{0}\otimes F_{1}+E_{1}\otimes
F_{0}.
\end{eqnarray*}

Check of $\ref{23}$ : Using again $K=e_{0}-e_{1}$, we have

\begin{eqnarray*}
S(e_{0}) &=&\frac{1}{2}S(\mathbf{1}+K)=\frac{1}{2}(\mathbf{1}+K^{-1})=\frac{1%
}{2}(\mathbf{1}+K^{-1})=e_{0} \\
S(e_{1}) &=&\frac{1}{2}S(\mathbf{1}-K)=\frac{1}{2}(\mathbf{1}-K^{-1})=\frac{1%
}{2}(\mathbf{1}-K^{-1})=e_{1}.
\end{eqnarray*}

Since $S$ is antiisomorphism, we also have

\begin{eqnarray*}
S(E_{0}) &=&S(e_{0}E)=S(E)S(e_{0})=(-K^{-1}E)e_{0}=-\left(
e_{0}-e_{1}\right) Ee_{0}=E_{1}, \\
S(E_{1}) &=&S(e_{1}E)=S(E)S(e_{1})=(-K^{-1}E)e_{1}=-\left(
e_{0}-e_{1}\right) Ee_{1}=-E_{0}, \\
S(F_{0}) &=&S(e_{0}F)=S(F)S(e_{0})=(-FK)e_{0}=-F\left( e_{0}-e_{1}\right)
e_{0}=-F_{1}, \\
S(F_{1}) &=&S(e_{1}F)=S(F)S(e_{1})=(-FK)e_{1}=-F\left( e_{0}-e_{1}\right)
e_{1}=F_{0}, \\
S(C_{0}) &=&S(e_{0}EF)=S(F)S(E)S(e_{0})=F\left( e_{0}-e_{1}\right) \left(
e_{0}-e_{1}\right) Ee_{0}=Fe_{1}E=e_{0}EF=C_{0}, \\
S(C_{1}) &=&S(e_{1}EF)=S(F)S(E)S(e_{1})=F\left( e_{0}-e_{1}\right) \left(
e_{0}-e_{1}\right) Ee_{1}=Fe_{0}E=e_{1}EF=C_{1}.
\end{eqnarray*}

Check of $\ref{24}$ : using again $K=e_{0}-e_{1}$, we have

\begin{eqnarray*}
\varepsilon (e_{0}) &=&\frac{1}{2}\varepsilon (\mathbf{1}+K)=\frac{1}{2}%
(1+1)=1, \\
\varepsilon (e_{1}) &=&\frac{1}{2}\varepsilon (\mathbf{1}-K)=\frac{1}{2}%
(1-1)=0.
\end{eqnarray*}

Since $\varepsilon $ is a morphism, we get, typically:,

$$
\varepsilon (E_{0})=\varepsilon (e_{0})\varepsilon (E)=1\times 0=0. 
$$
\edemo

\begin{remark}
The eigenvalues of $S$ in $H_{1}^{i}$ are $1$, $i$ and  $-i$ with the eigenspaces 
$$
V_{1}\;spanned\;by\;\left\{ 
\begin{array}{c}
e_{0} \\ 
e_{1} \\ 
C_{0} \\ 
C_{1}
\end{array}
\right. ,V_{i}\;spanned\;by\;\left\{ 
\begin{array}{c}
E_{0}-iE_{1} \\ 
F_{0}+iF_{1}
\end{array}
\right. and\;V_{-i}\;spanned\;by\;\left\{ 
\begin{array}{c}
E_{0}+iE_{1} \\ 
F_{0}-iF_{1}
\end{array}
.\right. 
$$
\end{remark}

\subsection{Hopf structure of $H_{2}^{i}$}

Recall that the algebra $H_{2}^{i}$ is defined by the symbols $K,E,F$ and the 
relations

$$
\left\{ 
\begin{array}{c}
KE=-EK, \\ 
KF=-FK, \\ 
\lbrack E,F]=\frac{1}{2i}(K-K^{-1}),
\end{array}
\right. \qquad \left\{ 
\begin{array}{c}
E^{2}=0, \\ 
F^{2}=0, \\ 
K^{4}=\mathbf{1},
\end{array}
\right. 
$$

or else symbols $E$, $F$ and $e_{m}=\frac{1}{4}\dsum\limits_{k\in \zzz/4\zzz}
i^{mk}K^{k}$, $i=0,1,2,3$, and relations
$$
\left\{ 
\begin{array}{c}
\dsum\limits_{m=0}^{3}e_{m}=\mathbf{1}, \\ 
e_{l}e_{m}=\delta _{lm},\quad l,m=0,1,2,3
\end{array}
\right. 
\left\{ 
\begin{array}{c}
E^{2}=0, \\ 
F^{2}=0, \\ 
\lbrack E,F]=e_{3}-e_{1},
\end{array}
\right. 
\left\{ 
\begin{array}{c}
e_{j}E=Ee_{j+2}, \\ 
Fe_{j}=e_{j+2}F,
\end{array}
\right.\quad j=0,1,2,3. 
$$
It is spanned by $e_{0},e_{2},E_{0},E_{2},F_{0},F_{2},P_{0},P_{2},e_{1},e_{3},E_{1},E_{3},F_{1},F_{3},P_{1},P_{3}$.

\begin{proposition} We have, for $m\in \zzz/4\zzz$,
\begin{equation}
\left\{ 
\begin{array}{l}
\Delta (e_{m})=\frac{1}{4}\dsum\limits_{k\in \zzz/4\zzz}e_{k}\otimes
e_{m-k}, \\ 
\Delta (E_{m})=\frac{1}{4}\dsum\limits_{k\in \zzz/4\zzz}\left[
E_{k}\otimes e_{m-k}+(i)^{-k}e_{k}\otimes E_{m-k}\right],  \\ 
\Delta (F_{m})=\frac{1}{4}\dsum\limits_{k\in \zzz/4\zzz}\left[
(i)^{m-k}F_{k}\otimes e_{m-k}+e_{k}\otimes F_{m-k}\right],  \\ 
\Delta (P_{m})=\frac{1}{4}\dsum\limits_{k\in \zzz/4\zzz}\left[
i^{m-k}P_{k}\otimes e_{m-k}+E_{k}\otimes F_{m-k}-(-1)^{k}(i)^{m}F_{k}\otimes
E_{m-k}+(i)^{-k}e_{k}\otimes P_{m-k}\right], 
\end{array}
\right.   \label{26}
\end{equation}
\begin{equation}
\left\{ 
\begin{array}{l}
S(e_{m})=e_{-m}, \\ 
S(E_{m})=(-i)^{m}E_{2-m}, \\ 
S(F_{m})=-(i)^{m}F_{2-m}, \\ 
S(P_{m})=P_{-m}-\delta _{3m}e_{1}+\delta _{1m}e_{3},
\end{array}
\right.   \label{27}
\end{equation}
\begin{equation}
\left\{
\begin{array}{l}
\varepsilon (e_{m}) =\delta _{0m}  \label{28} \\
\varepsilon (E_{m}) =\varepsilon (F_{m})=\varepsilon (P_{m})=0
\end{array}\right.  \nonumber
\end{equation}
\end{proposition}

\demo
Check of $\ref{26}$: Using $K=e_{0}-ie_{1}-e_{2}+ie_{3}$, $\dsum\limits_{k\in 
\zzz/4\zzz}(i)^{kl}=4\delta _{0l}$ for $l\in \zzz/4\zzz$, $%
Ke_{k}=e_{k}K=(i)^{-k}e_{k}$, $K^{-1}e_{k}=e_{k}=K^{-1}(i)^{k}e_{k}$, we have

\begin{eqnarray*}
\Delta (e_{m}) &=&\frac{1}{4}\dsum\limits_{r\in \zzz/4\zzz%
}(i)^{mr}\Delta (K)^{r}=\frac{1}{4}\dsum\limits_{l\in \zzz/4\zzz%
}(i)^{mr}K^{r}\otimes K^{r}=\frac{1}{4}\dsum\limits_{r,k,l\in \zzz/4\Bbb{Z%
}}(i)^{mr}\left( (i)^{-kr}e_{k}\otimes (i)^{-lr}e_{l}\right) \\
&=&\frac{1}{4}\dsum\limits_{r,k,l\in \zzz/4\zzz}(i)^{r(m-k-l)}\left(
e_{k}\otimes e_{l}\right) =\frac{1}{4}\dsum\limits_{k,l\in \zzz/4\zzz%
}\delta _{m,k+l}\left( e_{k}\otimes e_{l}\right) =\frac{1}{4}%
\dsum\limits_{k\in \zzz/4\zzz}\left( e_{k}\otimes e_{m-k}\right),
\end{eqnarray*}

\begin{eqnarray*}
\Delta (E_{m}) &=&\Delta (e_{m})\Delta (E)=\frac{1}{4}\left(
\dsum\limits_{k\in \zzz/4\zzz}e_{k}\otimes e_{m-k}\right) \left(
E\otimes \mathbf{1}+K\otimes E\right) \\
&=&\frac{1}{4}\dsum\limits_{k\in \zzz/4\zzz}\left( E_{k}\otimes
e_{m-k}+e_{k}K\otimes E_{m-k}\right) =\frac{1}{4}\dsum\limits_{k\in \zzz/4%
\zzz}\left( E_{k}\otimes e_{m-k}+(-i)^{k}e_{k}\otimes E_{m-k}\right),
\end{eqnarray*}

\begin{eqnarray*}
\Delta (F_{m}) &=&\Delta (e_{m})\Delta (F)=\frac{1}{4}\left(
\dsum\limits_{k\in \zzz/4\zzz}e_{k}\otimes e_{m-k}\right) \left(
F\otimes K^{-1}+\mathbf{1}\otimes F\right) \\
&=&\frac{1}{4}\dsum\limits_{k\in \zzz/4\zzz}\left( F_{k}\otimes
e_{m-k}K^{-1}+e_{k}\otimes F_{m-k}\right) =\frac{1}{4}\dsum\limits_{k\in 
\zzz/4\zzz}\left( (i)^{m-k}F_{k}\otimes e_{m-k}+e_{k}\otimes
F_{m-k}\right),
\end{eqnarray*}

\begin{eqnarray*}
\Delta (P_{m}) &=&\Delta (E_{m})\Delta (F)=\frac{1}{4}\dsum\limits_{k\in 
\zzz/4\zzz}\left( E_{k}\otimes e_{m-k}+(-i)^{k}e_{k}\otimes
E_{m-k}\right) \left( F\otimes K^{-1}+\mathbf{1}\otimes F\right) \\
&=&\frac{1}{4}\dsum\limits_{k\in \zzz/4\zzz}\left( E_{k}F\otimes
e_{m-k}K^{-1}+E_{k}\otimes F_{m-k}+(-i)^{k}F_{k}\otimes
E_{m-k}K^{-1}+(-i)^{k}e_{k}\otimes E_{m-k}F\right) \\
&=&\frac{1}{4}\dsum\limits_{k\in \zzz/4\zzz}\left(
(i)^{m-k}P_{k}\otimes e_{m-k}+E_{k}\otimes
F_{m-k}-(-1)^{k}(i)^{m}F_{k}\otimes E_{m-k}+(-i)^{k}e_{k}\otimes
P_{m-k}\right).
\end{eqnarray*}

Check of $\ref{27}$: We have by the antimultiplicativity of $S$,

$$
S(e_{m})=\dsum\limits_{k\in \zzz/4\zzz}(i)^{mk}S(K)^{k}=\dsum%
\limits_{k\in \zzz/4\zzz}(i)^{mk}K^{-k}=\dsum\limits_{k\in \zzz/4%
\zzz}(i)^{-mk}K^{k}=e_{-m}, 
$$

$$
S(E_{m})=S(E)S(e_{m})=(-K^{-1}E)e_{-m}=-K^{-1}e_{2-m}E=-(i)^{2-m}E_{2-m}=(-i)^{m}E_{2-m}, 
$$

$$
S(F_{m})=S(F)S(e_{m})=(-FK)e_{-m}=-(i)^{m}Fe_{-m}=-(i)^{m}e_{2-m}F=-(i)^{m}F_{2-m}, 
$$

\begin{eqnarray*}
S(P_{m})
&=&S(F)S(E)S(e_{m})=(-FK)(-K^{-1}E)e_{-m}=FEe_{-m}=(EF-e_{1}+e_{3})e_{-m} \\
&=&P_{-m}-\delta _{1,-m}e_{1}+\delta _{3,-m}e_{3}=P_{-m}-\delta
_{3,m}e_{1}+\delta _{1,m}e_{3}.
\end{eqnarray*}

Check of $\ref{28}$: Owing to the multiplicativity of $\varepsilon$, we have

$$
\varepsilon (e_{m})=\dsum\limits_{k\in \zzz/4\zzz}(i)^{mk}\varepsilon
(K)^{k}=\dsum\limits_{k\in \zzz/4\zzz}(i)^{mk}\delta _{0,m}, 
$$

whilst $\varepsilon (E_{m}),\varepsilon (F_{m})$ and $\varepsilon (P_{m})$
all vanish because they all contain the factor $\varepsilon (E)$ or $\varepsilon (F)$,
\edemo

\section{Adjoint representations of $H_{1}^{i}$ and $H_{2}^{i}$}


\subsection{Adjoint representation and adjoint trace of $H_{1}^{i}$}

\begin{proposition}
1) We have the following values of the biregular (or adjoint )
representation $\mu$
\begin{center}
\begin{tabular}{|c|c|c|c|c|c|c|c|c|}
\hline
$a$ & $\mu (a)e_{0}$ & $\mu (a)e_{1}$ & $\mu (a)E_{0}$ & $\mu (a)E_{1}$ & $\mu (a)F_{0}$
& $\mu (a)F_{1}$ & $\mu (a)C_{0}$ & $\mu (a)C_{1}$ \\ 
\hline
$e_{0}$ & $e_{0}$ & $e_{1}$ & 0 & 0 & 0 & 0 & $C_{0}$ & $C_{1}$ \\ 
\hline
$e_{1}$ & 0 & 0 & $E_{0}$ & $E_{1}$ & $F_{0}$ & $F_{1}$ & 0 & 0 \\ 
\hline
$E_{0}$ & 0 & 0 & 0 & 0 & $C_{0}+C_{1}$ & $C_{0}+C_{1}$ & 0 & 0 \\ 
\hline
$E_{1}$ & $-E_{0}+E_{1}$ & $E_{0}-E_{1}$ & 0 & 0 & 0 & 0 & 0 & 0 \\ 
\hline
$F_{0}$ & 0 & 0 & $-C_{0}-C_{1}$ & $C_{0}+C_{1}$ & 0 & 0 & 0 & 0 \\ 
\hline
$F_{1}$ & $F_{0}+F_{1}$ & $-F_{0}+F_{1}$ & 0 & 0 & 0 & 0 & 0 & 0 \\ 
\hline
$C_{0}$ & $2\left( C_{0}+C_{1}\right)$  & $2\left( C_{0}-C_{1}\right)$  & 0 & 0 & 0
& 0 & 0 & 0 \\ 
\hline
$C_{1}$ & 0 & 0 & 0 & 0 & 0 & 0 & 0 & 0\\
\hline
\end{tabular}
\end{center}

2) We have the following values of the adjoint trace $Tr_{\mu }$ of $H_{1}^{i}$
\begin{center}
\begin{tabular}{|c|c|c|c|c|c|c|c|c|}
\hline
$a$ & $e_{0}$ & $e_{1}$ & $E_{0}$ & $E_{1}$ & $F_{0}$ & $F_{1}$ & $C_{0}$ & $C_{1}$ \\ 
\hline
$tr_{\mu }(a)$ & 4 & 4 & 0 & 0 & 0 & 0 & 0 & 0\\
\hline
\end{tabular}
\end{center}

Thus, $Tr_{\mu}$ has the nilradical as its kernel and passes to the semi-simple quotient
as its trace.

3) The scalar product determined by the adjoint trace and the $*$-operation 
$$
<a,b>=\frac{1}{4}Tr_{\mu }(a^{*}b)\quad a,b\in H_{1}^{i}
$$

has the only non-vanishing values $<e_{0},e_{0}>=1$ and $<e_{1},e_{1}>=1$,
In other terms the scalar product $<.,.>$ is positive semi-definite with
null-space the nilradical $N_{1}^{i}$.
\end{proposition}

\demo
1) From $\ref{22}$, $\ref{23}$, and $\ref{24}$, we deduce the following table of
elements $(id\otimes S)\Delta (a)$ of $End_{\ccc}(H_{1}^{i})\otimes End_{%
\ccc}(H_{1}^{i})$, $a\in H_{1}^{i}$
\begin{eqnarray}
(id\otimes S)\Delta (e_{0}) &=&e_{0}\otimes e_{0}+e_{1}\otimes e_{1},
\label{30} \\
(id\otimes S)\Delta (e_{1}) &=&e_{0}\otimes e_{1}+e_{1}\otimes e_{0}, 
\nonumber \\
(id\otimes S)\Delta (E_{0}) &=&E_{0}\otimes e_{0}+E_{1}\otimes
e_{1}+e_{0}\otimes E_{1}+e_{1}\otimes E_{0},  \nonumber \\
(id\otimes S)\Delta (E_{1}) &=&E_{0}\otimes e_{1}+E_{1}\otimes
e_{0}-e_{0}\otimes E_{0}-e_{1}\otimes E_{1},  \nonumber \\
(id\otimes S)\Delta (F_{0}) &=&F_{0}\otimes e_{0}-F_{1}\otimes
e_{1}-e_{0}\otimes F_{1}+e_{1}\otimes F_{0},  \nonumber \\
(id\otimes S)\Delta (F_{1}) &=&-F_{0}\otimes e_{1}+F_{1}\otimes
e_{0}+e_{0}\otimes F_{0}-e_{1}\otimes F_{1},  \nonumber \\
(id\otimes S)\Delta (C_{0}) &=&C_{0}\otimes e_{0}-C_{1}\otimes
e_{1}+e_{0}\otimes C_{0}-e_{1}\otimes C_{1},  \nonumber \\
&&-E_{0}\otimes F_{1}+E_{1}\otimes F_{0}-F_{0}\otimes E_{1}+F_{1}\otimes
E_{0},  \nonumber \\
(id\otimes S)\Delta (C_{1}) &=&-C_{0}\otimes e_{1}+C_{1}\otimes
e_{0}+e_{0}\otimes C_{1}-e_{1}\otimes C_{0},  \nonumber \\
&&-F_{0}\otimes E_{0}+F_{1}\otimes E_{1}+E_{0}\otimes F_{0}-E_{1}\otimes
F_{1}.  \nonumber
\end{eqnarray}

From these relations, the corresponding $\mu (a)x,$ $x\in H_{1}^{i}$, are obtained by
making $\otimes \rightarrow x$ (observe that $\mu
(e_{0})$ and $\mu (e_{1})$ are the projections onto the even and odd
part of $End_{\ccc}(H_{1}^{i})$ for the $\zzz/2\zzz$-grading).

2) Results by inspection conferring $\ref{30}$ with $\ref{16}$.

3) The scalar product $<.,.>$ is given by the following table (where $<a,b>$
is plotted at the intersection of the line $a$ and the column $b$), obtained
by replacing the entries of \ref{16} by their quantum traces.

\begin{center}
\begin{tabular}{|c|c|c|c|c|c|c|c|c|}
\hline
& $e_{0}$ & $e_{1}$ & $E_{0}$ & $E_{1}$ & $F_{0}$ & $F_{1}$ & $C_{0}$ & $C_{1}$ \\ 
\hline
$e_{0}$ & 1 & 0 & 0 & 0 & 0 & 0 & 0 & 0 \\ 
\hline
$e_{1}$ & 0 & 1 & 0 & 0 & 0 & 0 & 0 & 0 \\ 
\hline
$E_{0}$ & 0 & 0 & 0 & 0 & 0 & 0 & 0 & 0 \\ 
\hline
$E_{1}$ & 0 & 0 & 0 & 0 & 0 & 0 & 0 & 0 \\ 
\hline
$F_{0}$ & 0 & 0 & 0 & 0 & 0 & 0 & 0 & 0 \\ 
\hline
$F_{1}$ & 0 & 0 & 0 & 0 & 0 & 0 & 0 & 0 \\ 
\hline
$C_{0}$ & 0 & 0 & 0 & 0 & 0 & 0 & 0 & 0 \\ 
\hline
$C_{1}$ & 0 & 0 & 0 & 0 & 0 & 0 & 0 & 0\\
\hline
\end{tabular}
\end{center}
\edemo

\subsection{Adjoint representation and adjoint trace of $H_{2}^{i}$}

\begin{proposition}
1) The biregular representation $\mu =\lambda *\rho $ of $H_{2}^{i}$
vanishes on $H_{2}^{i(1)}$ and is given on $H_{2}^{i(0)}$ by following table.
\begin{center}
\begin{tabular}{|c|c|c|c|c|c|c|c|c|}
\hline
$a$ & $\mu (a)e_{0}$ & $\mu (a)e_{2}$ & $\mu (a)E_{0}$ & $\mu (a)E_{2}$ & $\mu (a)F_{0}$
& $\mu (a)F_{2}$ & $\mu (a)P_{0}$ & $\mu (a)P_{2}$ \\ 
\hline
$_{0}$ & $\frac{1}{4}e_{0}$ & $\frac{1}{4}e_{2}$ & 0 & 0 & 0 & 0 & $\frac{1}{4}P_{0}$ & $\frac{1}{4}P_{2}$ \\ 
\hline
$e_{2}$ & 0 & 0 & $\frac{1}{4}E_{0}$ & $\frac{1}{4}E_{2}$ & $\frac{1}{4}F_{0}$ & 
$\frac{1}{4}F_{2}$ & 0 & 0 \\ 
\hline
$E_{0}$ & $\frac{1}{4}\left( E_{0}+E_{2}\right)$  & $\frac{1}{4}\left(
E_{0}+E_{2}\right)$  & 0 & 0 & 0 & 0 & 0 & 0 \\ 
\hline
$E_{2}$ & 0 & 0 & 0 & 0 & 0 & 0 & 0 & 0 \\ 
\hline
$F_{0}$ & 0 & 0 & $\frac{1}{2}P_{0}$ & $-\frac{1}{2}P_{2}$ & 0 & 0 & $\frac{1}{4}P_{0}$ & $-\frac{1}{4}P_{2}$ \\ 
\hline
$F_{2}$ & $\frac{1}{4}\left( F_{0}+F_{2}\right)$  & $\frac{1}{4}\left(
F_{0}+F_{2}\right)$  & 0 & 0 & 0 & 0 & 0 & 0 \\ 
\hline
$P_{0}$ & 0 & 0 & 0 & 0 & 0 & 0 & 0 & 0 \\ 
\hline
$P_{2}$ & 0 & 0 & 0 & 0 & 0 & 0 & 0 & 0\\
\hline
\end{tabular}
\end{center}

and

\begin{center}
\begin{tabular}{|c|c|c|c|c|c|c|c|c|}
\hline
$a $& $\mu (a)e_{1} $& $\mu (a)e_{3}$ & $\mu (a)E_{1}$ & $\mu (a)E_{3}$ & $\mu (a)F_{1}$
& $\mu (a)F_{3}$ & $\mu (a)P_{1}$ & $\mu (a)P_{3}$ \\ 
\hline
$e_{0}$ & $\frac{1}{4}e_{1}$ & $\frac{1}{4}e_{3}$ & 0 & 0 & 0 & 0 & $\frac{1}{4}P_{1}$ & $\frac{1}{4}P_{3}$ \\ 
\hline
$e_{2}$ & 0 & 0 & $\frac{1}{4}E_{1}$ & $\frac{1}{4}E_{3}$ & $\frac{1}{4}F_{1}$ & $\frac{1}{4}F_{3}$ & 0 & 0 \\
\hline 
$E_{0}$ & $\frac{1}{4}\left( E_{1}+E_{3}\right)$  & $\frac{1}{4}\left(
E_{1}+E_{3}\right)$  & 0 & 0 & $\frac{1}{4}e_{1}$ & $-\frac{1}{4}e_{3}$ & $-\frac{1}{4}E_{1}$ & $\frac{1}{4}E_{3}$ \\
\hline 
$E_{2}$ & 0 & 0 & 0 & 0 & 0 & 0 & 0 & 0 \\ 
\hline
$F_{0}$ & 0 & 0 & $ie_{1}$ & $ie_{1}$ & 0 & 0 & $\frac{i}{4}P_{1}$ & $-\frac{i}{4}P_{3}$ \\ 
\hline
$F_{2}$ & $-\frac{i}{4}\left( F_{1}-F_{3}\right)$  & $\frac{i}{4}\left(
F_{1}-F_{3}\right)$  & 0 & 0 & 0 & 0 & 0 & 0 \\ 
\hline
$P_{0}$ & $\frac{1}{4}\left( e_{1}+e_{3}\right)$  & $\frac{1}{4}\left(
e_{1}+e_{3}\right) $ & 0 & 0 & 0 & 0 & 0 & 0 \\ 
\hline
$P_{2}$ & 0 & 0 & 0 & 0 & 0 & 0 & 0 & 0\\
\hline
\end{tabular}
\end{center}

We have the generic formulae 
\begin{equation}
\left\{ 
\begin{array}{c}
\mu (e_{m})e_{j}=\frac{1}{4}\delta _{0,m}e_{j}, \\ 
\mu (e_{m})E_{j}=\frac{1}{4}\delta _{0,m}E_{j}, \\ 
\mu (e_{m})F_{j}=\frac{1}{4}\delta _{2,m}F_{j}, \\ 
\mu (e_{m})P_{j}=\frac{1}{4}\delta _{0,m}P_{j},
\end{array}
\right. \left\{ 
\begin{array}{c}
\mu (E_{m})e_{j}=\frac{1}{4}\delta _{0,m}\left( E_{j+2}+E_{j}\right),  \\ 
\mu (E_{m})E_{j}=0, \\ 
\mu (E_{m})F_{j}=\frac{1}{4}\delta _{0,m}\left( \delta _{1,j}e_{1}-\delta
_{3,j}e_{3}\right),  \\ 
\mu (E_{m})P_{j}=\frac{1}{4}\delta _{0,m}\left( \delta _{1,j}E_{1}-\delta
_{3,j}E_{3}\right), 
\end{array}
\right.   \label{31}
\end{equation}
$$
\left\{ 
\begin{array}{c}
\mu (F_{m})e_{j}=\frac{1}{4}\delta _{2,m}\left( \left( i\right)
^{j}F_{j+2}+\left( i\right) ^{-j}F_{j}\right),  \\ 
\mu (F_{m})E_{j}=\frac{1}{4}\delta _{0,m}\left( \left( \left( i\right)
^{j}+\left( i\right) ^{-j}\right) P_{j}+i\left( \delta _{1,j}e_{1}+\delta
_{3,j}e_{3}\right) \right),  \\ 
\mu (F_{m})F_{j}=0, \\ 
\mu (F_{m})P_{j}=\frac{1}{4}\delta _{0,m}\left( i\right) ^{j}P_{j},
\end{array}
\right. 
$$

$$
\left\{ 
\begin{array}{c}
\mu (P_{m})e_{j}=\frac{1}{4}i\delta _{0,m}\left( \delta _{1,j}+\delta
_{3,j}\right) \left( e_{1}+e_{3}\right),  \\ 
\mu (P_{m})E_{j}=-\frac{1}{4}\delta _{0,m}\left( \delta _{1,j}E_{1}-\delta
_{3,j}E_{3}\right),  \\ 
\mu (P_{m})F_{j}=0, \\ 
\mu (P_{m})P_{j}=0.
\end{array}
\right. 
$$

2) We have the following values of the adjoint trace $Tr_{\mu }$ of $H_{2}^{i}$.

\begin{center}
\begin{tabular}{|c|c|c|c|c|c|c|c|c|c|c|c|c|c|c|c|c|}
\hline
$a$ & $e_{0}$ & $e_{2}$ & $E_{0} $& $E_{2}$ & $F_{0}$ & $F_{2}$ & $C_{0}$ & $C_{2}$ & $e_{1}$ & 
$e_{3}$ & $E_{1}$ & $E_{3}$ & $F_{1}$ & $F_{3}$ & $P_{1}$ & $P_{3}$ \\ 
\hline
$Tr_{\mu }(a)$ & 1 & 2 & 0 & 0 & 0 & 0 & 0 & 0 & 0 & 0 & 0 & 0 & 0 & 0 & 0 & 0\\
\hline
\end{tabular}
\end{center}

3) The scalar product determined by the trace and the $*$-operation 
$$
<a,b>=\frac{1}{4}Tr_{\mu }(a^{*}b)\quad a,b\in H_{2}^{i}
$$

has the only non vanishing values $<e_{0},e_{0}>=1$ and $<e_{2},e_{2}>=2$. In other terms, the scalar product $<.,.>$ is positive semi definite with
null-space $N_{2}^{i(0)}\oplus H_{2}^{i(1)}$.
\end{proposition}

\demo
1) From $\ref{26}$ and $\ref{27}$ we deduce the following elements $\left(
id\otimes S\right) \Delta (a)$ of $End_{\ccc}(H_{2}^{i})\otimes End_{\Bbb{%
C}}(H_{2}^{i})$, $a\in H_{2}^{i}$, $\left( id\otimes S\right) \Delta (e_{m})=%
\frac{1}{4}\dsum\limits_{k\in \zzz/2\zzz}e_{k}\otimes e_{k-m}$, from
which the corresponding $4\mu (a)x,x\in H_{2}^{i}$, are obtained by making $%
\otimes \rightarrow x$, yielding \ref{31}

$$
\left\{ 
\begin{array}{c}
\mu (e_{m})e_{j}=\frac{1}{4}\dsum\limits_{k\in \zzz/2\zzz%
}e_{k}e_{j}e_{k-m}=\frac{1}{4}\delta _{0,m}\dsum\limits_{k\in \zzz/2\Bbb{Z%
}}\delta _{k,j}e_{j}=\frac{1}{4}\delta _{0,m}e_{j}, \\ 
\mu (e_{m})E_{j}=\frac{1}{4}\dsum\limits_{k\in \zzz/2\zzz%
}e_{k}E_{j}e_{k-m}=\frac{1}{4}\delta _{0,m}\dsum\limits_{k\in \zzz/2\Bbb{Z%
}}\delta _{k,j}E_{j}=\frac{1}{4}\delta _{0,m}E_{j}, \\ 
\mu (e_{m})F_{j}=\frac{1}{4}\dsum\limits_{k\in \zzz/2\zzz%
}e_{k}F_{j}e_{k-m}=\frac{1}{4}\delta _{k+2,k-m}\dsum\limits_{k\in \zzz/2%
\zzz}\delta _{k,j}F_{j}=\frac{1}{4}\delta _{2,m}F_{j}, \\ 
\mu (e_{m})P_{j}=\frac{1}{4}\dsum\limits_{k\in \zzz/2\zzz%
}e_{k}P_{j}e_{k-m}=\frac{1}{4}\delta _{0,m}\dsum\limits_{k\in \zzz/2\Bbb{Z%
}}\delta _{k,j}p_{j}=\frac{1}{4}\delta _{0,m}P_{j}.
\end{array}
\right. 
$$

The derivation of the other formulae is obtained via multiplicativity of $%
\mu $. Combining the latter with

\begin{equation}
\left\{ 
\begin{array}{c}
\mu (E)(e_{j})=E_{j+2}+E_{j}, \\ 
\mu (E)(E_{j})=0, \\ 
\mu (E)(F_{j})=-\left( \delta _{1,j}e_{1}-\delta _{3,j}e_{3}\right), \\ 
\mu (E)(P_{j})=-\left( \delta _{1,j}E_{1}-\delta _{3,j}E_{3}\right),
\end{array}
\right. \left\{ 
\begin{array}{c}
\mu (F)(e_{j})=\left( i\right) ^{j}F_{j+2}+\left( i\right) ^{-j}F_{j}, \\ 
\mu (F)(E_{j})=-\left( \left( i\right) ^{j}+\left( i\right) ^{-j}\right)
P_{j}-\left( i\right) ^{j}\left( \delta _{1,j}e_{1}-\delta _{3,j}e_{3}\right),
\\ 
\mu (F)(F_{j})=0, \\ 
\mu (F)(P_{j})=-i\left( E_{1}+E_{3}\right),
\end{array}
\right.  \label{32}
\end{equation}

\begin{equation}
\left\{ 
\begin{array}{c}
\mu (P)(e_{j})=i\left( \delta _{1,j}+\delta _{3,j}\right) \left(
e_{1}+e_{3}\right), \\ 
\mu (P)(E_{j})=2i\left( E_{1}+E_{3}\right), \\ 
\mu (P)(F_{j})=0, \\ 
\mu (P)(P_{j})=0.
\end{array}
\right.  \label{33}
\end{equation}

they are checked as follows: from $\ref{2}$ we have $\left( id\otimes S\right)
\Delta (E)\mathbf{=}E\otimes \mathbf{1}+K\otimes K^{-1}E$ and $\left(
id\otimes S\right) \Delta (F)\mathbf{=}F\otimes K^{-1}-\mathbf{1}\otimes FK$.  Hence we havefor $x\in H_{2}^{i}$

\begin{eqnarray*}
\mu (E)(x) &=&Ex+KxK^{-1}E, \\
\mu (F)(x) &=&FxK^{-1}+xKF,
\end{eqnarray*}
and thus
\begin{eqnarray*}
\mu (E)e_{j} &=&Ee_{j}+Ke_{j}K^{-1}E=E_{j+2}+E_{j}, \\
\mu (E)E_{j} &=&EE_{j}+KE_{j}K^{-1}E=0, \\
\mu (E)F_{j} &=&EF_{j}+KF_{j}K^{-1}E=P_{j}-e_{j}\left( EF+e_{1}-e_{3}\right)
=-\left( \delta _{1,j}e_{1}-\delta _{3,j}e_{3}\right), \\
\mu (E)P_{j} &=&EP_{j}+KP_{j}K^{-1}E=e_{j}EFE=-e_{j}\left(
e_{1}-e_{3}\right) E=-\left( \delta _{1,j}E_{1}-\delta _{3,j}E_{3}\right),
\end{eqnarray*}
and

\begin{eqnarray*}
\mu (F)e_{j} &=&Fe_{j}K^{-1}+e_{j}KF=\left( i\right) ^{j}F_{j+2}+\left(
i\right) ^{-j}F_{j}, \\
\mu (F)E_{j} &=&FE_{j}K^{-1}+E_{j}KF=\left( i\right) ^{j+2}e_{j}\left(
EF+e_{1}-e_{3}\right) -\left( i\right) ^{-j}P_{j}, \\
&=&-\left( \left( i\right) ^{j}+\left( i\right) ^{-j}\right) P_{j}-\left(
i\right) ^{j}\left( \delta _{1,j}e_{1}-\delta _{3,j}e_{3}\right), \\
\mu (F)F_{j} &=&FF_{j}K^{-1}+F_{j}KF=0, \\
\mu (F)P_{j} &=&FP_{j}K^{-1}+P_{j}KF=\left( i\right) ^{j}e_{j+2}FEF=\left(
i\right) ^{j}e_{j+2}\left( e_{1}-e_{3}\right) F, \\
&=&\left( i\right) ^{j}\left( \delta _{3,j}F_{1}-\delta _{1,j}F_{3}\right)
=-i\left( E_{1}+E_{3}\right),
\end{eqnarray*}
further
\begin{eqnarray*}
\mu (P)e_{j} &=&\mu (E)\mu (F)e_{j}=\mu (E)\left( \left( i\right)
^{j}F_{j+2}+\left( i\right) ^{-j}F_{j}\right), \\
&=&-\left( i\right) ^{j}\left( \delta _{1,j+2}e_{1}-\delta
_{3,j+2}e_{3}\right) -\left( i\right) ^{-j}\left( \delta _{1,j}e_{1}-\delta
_{3,j}e_{3}\right), \\
&=&-\left( i\right) ^{j}\left( \delta _{3,j}e_{1}-\delta _{1,j}e_{3}\right)
-\left( i\right) ^{-j}\left( \delta _{1,j}e_{1}-\delta _{3,j}e_{3}\right), \\
&=&i\delta _{3,j}e_{1}+i\delta _{1,j}e_{3}+i\delta _{1,j}e_{1}+i\delta
_{3,j}e_{3}=i\left( \delta _{1,j}+\delta _{3,j}\right) \left(
e_{1}+e_{3}\right), \\
\mu (P)E_{j} &=&\mu (E)\mu (F)E_{j}=\mu (E)\left( -\left( \left( i\right)
^{j}+\left( i\right) ^{-j}\right) P_{j}-\left( i\right) ^{j}\left( \delta
_{1,j}e_{1}-\delta _{3,j}e_{3}\right) \right), \\
&=&\left( \left( i\right) ^{j}+\left( i\right) ^{-j}\right) \left( \delta
_{1,j}E_{1}-\delta _{3,j}E_{3}\right) +i\mu (E)\left( e_{1}+e_{3}\right), \\
&=&2i\left( E_{1}+E_{3}\right). \\
\mu (p)F_{j} &=&\mu (E)\mu (F)F_{j}=0
\end{eqnarray*}

We now compute successively

\begin{eqnarray*}
\mu (E_{m})e_{j} &=&\mu (e_{m})\mu (E)e_{j}=\mu (e_{m})\left(
E_{j+2}+E_{j}\right) =\frac{1}{4}\delta _{0,m}\left( E_{j+2}+E_{j}\right), \\
\mu (E_{m})E_{j} &=&\mu (E)\mu (e_{m+2})E_{j}=\frac{1}{4}\delta _{0,m+2}\mu
(E)(E_{j})=0, \\
\mu (E_{m})F_{j} &=&\mu (E)\mu (e_{m+2})F_{j}=\frac{1}{4}\delta _{2,m+2}\mu
(E)(F_{j})=\frac{1}{4}\delta _{0,m}\left( \delta _{1,j}e_{1}-\delta
_{3,j}e_{3}\right), \\
\mu (E_{m})P_{j} &=&\mu (e_{m})\mu (E)P_{j}=\mu (E), \\
&=&-\mu (e_{m})\left( \delta _{1,j}E_{1}-\delta _{3,j}E_{3}\right) =-\frac{1%
}{4}\delta _{0,m}\left( \delta _{1,j}E_{1}-\delta _{3,j}E_{3}\right),
\end{eqnarray*}

and

\begin{eqnarray*}
\mu (F_{m})e_{j} &=&\mu (e_{m})\mu (F)e_{j}=\mu (e_{m})\left( \left(
i\right) ^{j}F_{j+2}+\left( i\right) ^{-j}F_{j}\right) =\frac{1}{4}\delta
_{2,m}\left( \left( i\right) ^{j}F_{j+2}+\left( i\right) ^{-j}F_{j}\right), \\
\mu (F_{m})E_{j} &=&\mu (e_{m})\mu (F)E_{j}=\mu (e_{m})\left( -\left( \left(
i\right) ^{j}+\left( i\right) ^{-j}\right) P_{j}-\left( i\right) ^{j}\left(
\delta _{1,j}e_{1}-\delta _{3,j}e_{3}\right) \right), \\
&=&-\frac{1}{4}\delta _{0,m}\left( \left( \left( i\right) ^{j}+\left(
i\right) ^{-j}\right) P_{j}+\left( i\right) ^{j}\left( \delta
_{1,j}e_{1}-\delta _{3,j}e_{3}\right) \right), \\
\mu (F_{m})F_{j} &=&\mu (e_{m})\mu (F)F_{j}=0, \\
\mu (F_{m})P_{j} &=&\mu (e_{m})\mu (F)P_{j}=\left( i\right) ^{j}\mu
(e_{m})P_{j}=\frac{1}{4}\delta _{0,m}\left( i\right) ^{j}P_{j},
\end{eqnarray*}

and

\begin{eqnarray*}
\mu (P_{m})e_{j} &=&\mu (e_{m})\mu (P)e_{j}=i\mu (e_{m})\left( \delta
_{1,j}+\delta _{3,j}\right) \left( e_{1}+e_{3}\right), \\
&=&\frac{1}{4}i\delta _{0,m}\left( \delta _{1,j}+\delta _{3,j}\right) \left(
e_{1}+e_{3}\right), \\
\mu (P_{m})E_{j} &=&\mu (e_{m})\mu (P)E_{j}=-\mu (e_{m})\left( \delta
_{1,j}E_{1}-\delta _{3,j}E_{3}\right) =-\frac{1}{4}\delta _{0,m}\left(
\delta _{1,j}E_{1}-\delta _{3,j}E_{3}\right), \\
\mu (P_{m})F_{j} &=&\mu (E)\mu (e_{m+2})F_{j}=0, \\
\mu (P_{m})P_{j} &=&\mu (e_{m})\mu (E)\mu (F)P_{j}=0.
\end{eqnarray*}
\edemo

\section{Idempotents, automorphisms and real forms of $H_{1}^{i}$}

\begin{lemma}
1) There is a unique $\kappa \in Aut(H_{1}^{i})$ (the flip) such that 
$$
\left\{ 
\begin{array}{c}
\kappa (K)=-K, \\ 
\kappa (E)=E, \\ 
\kappa (F)=F.
\end{array}
\right. 
$$

2) $\kappa $ is an involution performing the exchanges $e_{0}\leftrightarrow
e_{1},E_{0}\leftrightarrow E_{1},F_{0}\leftrightarrow F_{1}$ and $%
C_{0}\leftrightarrow C_{1}$.

3) One has $\kappa \circ S\circ \kappa =S^{-1}$.
\end{lemma}

\demo
1) and 2) : immediate from $\ref{16}$.

3) straightforward from $S^{-1}(K)=K^{-1}$, $S^{-1}(E)=K^{-1}E$ and $S^{-1}(K)=FK$.
\edemo

\subsection{Idempotents of $H_{1}^{i}$}

\begin{proposition}
1) The idempotents of $H_{1}^{i}$ are

-the element $0$ of rank $0$,

-the unit $\mathbf{1}$ of rank $8$,

-a continous famuly of rank 4 idempotents 
$$
\left\{ 
\begin{array}{c}
e_{0,\beta ,\gamma ,\delta ,\eta }=e_{0}+\beta E_{0}+\gamma F_{0}+\delta
E_{1}+\eta F_{1}+\left( \beta \eta +\delta \gamma \right) \left(
C_{1}-C_{0}\right),  \\ 
e_{1,\beta ,\gamma ,\delta ,\eta }=e_{1}+\beta E_{0}+\gamma F_{0}+\delta
E_{1}+\eta F_{1}+\left( \beta \eta +\delta \gamma \right) \left(
C_{1}-C_{0}\right), 
\end{array}
\right. 
$$

indexed by the parameters $\beta ,\gamma ,\delta ,\eta \in \ccc$.

2) Consequently the automorphisms $\phi $ of $H_{1}^{i}$ are of either the two types:   

a)fulfilling $\phi (e_{0})=e_{0,\beta ,\gamma ,\delta ,\eta },\phi
(e_{1})=e_{1,-\beta ,-\gamma ,-\delta ,-\eta }$,

b)fulfilling $\phi (e_{0})=e_{1,\beta ,\gamma ,\delta ,\eta },\phi
(e_{1})=e_{0,-\beta ,-\gamma ,-\delta ,-\eta }$.
\end{proposition}

\demo
1) Every idempotent $e$ of $H_{1}^{i}=F_{1}^{i}\oplus N_{1}^{i}$ ($F_{1}^{i}$
a subalgebra, $N_{1}^{i}$ an ideal) decomposes as $e=e^{\prime }+e^{\prime
\prime }$ with $e^{\prime 2}=e^{\prime },e^{\prime }e^{\prime \prime
}+e^{\prime \prime }e^{\prime }+e^{\prime \prime 2}=e^{\prime \prime }$. We
search $e^{\prime }$

$e^{\prime }=\alpha e_{0}+\beta e_{1}=e^{\prime 2}=\alpha ^{2}e_{0}+\beta
^{2}e_{1}$ yields $\alpha \left( \alpha -1\right) =\beta \left( \beta
-1\right) =0$, hence $e^{\prime }$ either $0,e_{0},e_{1},$ or $e_{0}+e_{1}$.

We search $e^{\prime \prime
}=x_{0}E_{0}+y_{0}F_{0}+v_{0}C_{0}+x_{1}E_{1}+y_{1}F_{1}+v_{1}C_{1}$ with $%
e^{\prime \prime }=(x_{0}y_{1}+x_{1}y_{0})(C_{0}+C_{1})$

-To $e^{\prime }=0,e^{\prime \prime 2}=e^{\prime \prime }$ yields $%
x_{0}=y_{0}=x_{1}=y_{1}=0,v_{0}=v_{1}=x_{0}y_{1}+x_{1}y_{0}=0$ thus $%
e^{\prime \prime }=0$

-To $e^{\prime }=\mathbf{1},2e^{\prime \prime }+e^{\prime \prime
2}=e^{\prime \prime },e^{\prime \prime 2}=-e^{\prime \prime }$ yields $%
x_{0}=y_{0}=x_{1}=y_{1}=0$, $v_{0}=v_{1}=-\left(
x_{0}y_{1}+x_{1}y_{0}\right) x_{0}y_{1}+x_{1}y_{0}=0$ thus $e^{\prime \prime
}=0$

-To $e^{\prime }=e_{0}$ since $e_{0}e^{\prime \prime }+e_{0}e^{\prime \prime
}=x_{0}E_{0}+y_{0}F_{0}+v_{0}C_{0}+x_{1}E_{1}+y_{1}F_{1}+v_{1}C_{1}=e^{%
\prime \prime }++v_{0}C_{0}-v_{1}C_{1}$

$e_{0}e^{\prime \prime }+e_{0}e^{\prime \prime }=e^{\prime \prime }$ yields $%
v_{1}=-v_{0}=x_{0}y_{1}+x_{1}y_{0}$, $x_{0},y_{0},x_{1},y_{1}$ remaining
arbitrary whence

$e=e_{0}+x_{0}E_{0}+y_{0}F_{0}+x_{1}E_{1}+y_{1}F_{1}+\left(
x_{0}y_{1}+x_{1}y_{0}\right) (C_{0}-C_{1})$. The flip then yields by
exchange symmetry $0\leftrightarrow 1$ the idempotent $%
e_{1}+x_{0}E_{0}+y_{0}F_{0}+x_{1}E_{1}+y_{1}F_{1}+\left(
x_{0}y_{1}+x_{1}y_{0}\right) (C_{0}-C_{1})$

Rank of $e_{0,\beta ,\gamma ,\delta ,\eta }$ : $%
h=a_{0}e_{0}+X_{0}E_{0}+Y_{0}F_{0}+c_{0}C_{0}+a_{1}e_{1}+X_{1}E_{1}+Y_{1}F_{1}+c_{1}C_{1} 
$ fulfills $e_{0,\beta ,\gamma ,\delta ,\eta }h=h$ iff one has the relations 
$a_{1}=0$, $\eta X_{0}+\delta Y_{0}=(\beta \eta +\delta \gamma )a_{0}$
(automatic), $X_{1}=\delta a_{0}$, $Y_{1}=\eta a_{0}$. Thus in the two
occurring cases $\delta \eta =0$ and $\delta \neq 0,\eta \neq 0$, one has
rank 4.

2) Obvious from (1) since automorphisms turn an idempotent into an
idempotent of the same rank, Morover, $\varphi $ must be bijective and satisfy $\varphi
(e_{1})+\varphi (e_{0})=\mathbf{1}$.
\edemo

\subsection{Automorphisms of $H_{1}^{i}$}

\begin{proposition}
1) The set $Aut(H_{1}^{i})$ of automorphisms of $H_{1}^{i}$ coinciding with $%
Aut(N_{1}^{i})$ consists of elements of the two types (whose respective sets
will be denoted $Aut_{I}(H_{1}^{i})$ and $Aut_{II}(H_{1}^{i})$.

Type $I$ : 
\begin{equation}
\begin{array}[t]{l}
\varphi (e_{0})=e_{0}+\beta E_{0}+\gamma F_{0}+\delta E_{1}+\eta
F_{1}+\left( \beta \eta +\gamma \delta \right) \left( C_{1}-C_{0}\right),  \\ 
\varphi (E_{0})=\mu _{0}E_{0}+\nu _{0}F_{0}+\left( \delta \nu _{0}+\eta \mu
_{0}\right) \left( C_{1}-C_{0}\right),  \\ 
\varphi (F_{0})=\sigma _{0}E_{0}+\tau _{0}F_{0}+\left( \delta \tau _{0}+\eta
\sigma _{0}\right) \left( C_{1}-C_{0}\right),  \\ 
\varphi (C_{0})=\lambda C_{0},
\end{array}
\label{34}
\end{equation}

and

\begin{equation}
\begin{array}[t]{l}
\varphi (e_{1})=e_{1}-\beta E_{0}-\gamma F_{0}-\delta E_{1}-\eta
F_{1}-\left( \beta \eta +\gamma \delta \right) \left( C_{1}-C_{0}\right),  \\ 
\varphi (E_{1})=\mu _{1}E_{1}+\nu _{1}F_{1}+\left( \beta \nu _{1}+\gamma \mu
_{1}\right) \left( C_{1}-C_{0}\right),  \\ 
\varphi (F_{1})=\sigma _{1}E_{1}+\tau _{1}F_{1}+\left( \beta \tau
_{1}+\gamma \sigma _{1}\right) \left( C_{1}-C_{0}\right),  \\ 
\varphi (C_{1})=\lambda C_{1},
\end{array}
\label{35}
\end{equation}

for constants $\beta ,\gamma ,\delta ,\eta ,\lambda ,\mu _{0},\nu
_{0},\sigma _{0},\tau _{0},\mu _{1},\nu _{1},\sigma _{1},\tau _{1}\in \ccc
$ constrained by

\begin{equation}
\left\{ 
\begin{array}{l}
1)\mu _{0}\nu _{1}+\nu _{0}\mu _{1}=0, \\ 
2)\sigma _{0}\tau _{1}+\tau _{0}\sigma _{1}=0, \\ 
3)\lambda =\mu _{0}\tau _{1}+\nu _{0}\sigma _{1}=\sigma _{0}\nu _{1}+\tau
_{0}\mu _{1},
\end{array}
\right. \left\{ 
\begin{array}{l}
4)\mu _{0}\tau _{0}-\nu _{0}\sigma _{0}\neq 0, \\ 
5)\mu _{1}\tau _{1}-\nu _{0}\sigma _{0}\neq 0, \\ 
6)\lambda \neq 0.
\end{array}
\right.   \label{36}
\end{equation}

Type $II$ : Product $\varphi \kappa $ (or for that matter $\kappa \varphi $)
with $\varphi $ of the preceding type $I$%
$$
\begin{array}[t]{l}
\varphi (e_{0})=e_{1}+\beta E_{1}+\gamma F_{1}+\delta E_{0}+\eta
F_{0}+\left( \beta \eta +\gamma \delta \right) \left( C_{0}-C_{1}\right),  \\ 
\varphi (E_{0})=\mu _{0}E_{1}+\nu _{0}F_{1}+\left( \delta \nu _{0}+\eta \mu
_{0}\right) \left( C_{0}-C_{1}\right),  \\ 
\varphi (F_{0})=\sigma _{0}E_{1}+\tau _{0}F_{1}+\left( \delta \tau _{0}+\eta
\sigma _{0}\right) \left( C_{0}-C_{1}\right),  \\ 
\varphi (C_{0})=\lambda C_{1},
\end{array}
$$
and

$$
\begin{array}[t]{l}
\varphi (e_{1})=e_{0}-\beta E_{1}-\gamma F_{1}-\delta E_{0}-\eta
F_{0}-\left( \beta \eta +\gamma \delta \right) \left( C_{0}-C_{1}\right),  \\ 
\varphi (E_{1})=\mu _{1}E_{0}+\nu _{1}F_{0}+\left( \beta \nu _{1}+\gamma \mu
_{1}\right) \left( C_{0}-C_{1}\right),  \\ 
\varphi (F_{1})=\sigma _{1}E_{0}+\tau _{1}F_{0}+\left( \beta \tau
_{1}+\gamma \sigma _{1}\right) \left( C_{0}-C_{1}\right),  \\ 
\varphi (C_{1})=\lambda C_{0},
\end{array}
$$

with constants constrained as in $\ref{36}$. Observe that these constraints,
as well as the whole structure of $H_{1}^{i}$, is invariant under
the flip.

2) $Aut_{I}(H_{1}^{i})$ is a normal subgroup of $Aut(H_{1}^{i})$.

3) The subgroup $Aut(H_{1}^{i})=Aut_{I}(H_{1}^{i})\oplus Aut_{II}(H_{1}^{i})$
of $Aut_{I}(H_{1}^{i})$ specified by $\beta =\gamma =\delta =\eta =0$,
consists of the $4\times 4$ matrices 
$$
M=\left( 
\begin{array}{llll}
\mu _{0} & \nu _{0} & 0 & 0 \\ 
\sigma _{0} & \tau _{0} & 0 & 0 \\ 
0 & 0 & \mu _{1} & \nu _{1} \\ 
0 & 0 & \sigma _{1} & \tau _{1}
\end{array}
\right) \in M_{2}(\ccc)\oplus M_{2}(\ccc),
$$

leaving stable the set of bilinear forms with vanishing sum of elements of
their second diagonal. Specifically, for each matrix of the form

$$
G=\left( 
\begin{array}{llll}
0 & 0 & a & b \\ 
0 & 0 & -b & d \\ 
a & b & 0 & 0 \\ 
b & d & 0 & 0
\end{array}
\right) \quad a,b,c,d\in \ccc,
$$

$M^{t}GM$ is a matrix of the same type.

4) The constraints $\ref{36}$ entail the relations 
\begin{equation}
\left\{ 
\begin{array}{c}
\nu _{0}\tau _{1}+\tau _{0}\nu _{1}=0, \\ 
\mu _{0}\sigma _{1}+\sigma _{0}\mu _{1}=0.
\end{array}
\right.   \label{37}
\end{equation}

Note that (2) answers the natural question why the constraints $\ref{36}$
propagate through matrix products. Indeed it is not clear a priori that $%
Aut_{I}(H_{1}^{i})$ characterized as in (1) is multiplicative. Note also how
two symmetries permute (either respecting or exchanging) formula, the flip $%
(0\leftrightarrow 1)$ and the symmetry $(E\leftrightarrow F,\mu
\leftrightarrow \tau ,\nu \leftrightarrow \sigma )$.
\end{proposition}

\demo
1) Exhausting the constraints stemming from multiplicativity of $\varphi$,
in view of the multiplication table $\ref{16}$, we have the implications

\begin{eqnarray*}
\varphi (e_{0}E_{0}) &=&\varphi (e_{0})\varphi (E_{0})\Rightarrow \mu
_{1}^{\prime }=\nu _{1}^{\prime }=0,\rho _{1}^{\prime }=\delta \nu _{0}+\eta
\mu _{0}, \\
\varphi (e_{0}F_{0}) &=&\varphi (e_{0})\varphi (F_{0})\Rightarrow \sigma
_{1}^{\prime }=\tau _{1}^{\prime }=0,\omega _{1}^{\prime }=\delta \tau
_{0}+\eta \sigma _{0}, \\
\varphi (E_{0}e_{1}) &=&\varphi (E_{0})\varphi (e_{1})\Rightarrow \mu
_{1}^{\prime }=\nu _{1}^{\prime }=0,\rho _{0}^{^{\prime }}=-\left( \delta
\nu _{0}+\eta \mu _{0}\right), \\
\varphi (F_{0}e_{1}) &=&\varphi (F_{0})\varphi (e_{1})\Rightarrow \sigma
_{1}^{\prime }=\tau _{1}^{\prime }=0,\omega _{1}^{\prime }=-\left( \delta
\tau _{0}+\eta \sigma _{0}\right), \\
\varphi (e_{1}F_{1}) &=&\varphi (e_{1})\varphi (F_{1})\Rightarrow \sigma
_{0}^{^{\prime }}=\tau _{0}^{\prime }=0,\omega _{0}^{\prime }=-\left( \beta
\tau _{1}+\gamma \sigma _{1}\right), \\
\varphi (E_{1}e_{0}) &=&\varphi (E_{1})\varphi (e_{0})\Rightarrow \mu
_{0}^{\prime }=\nu _{0}^{\prime }=0,\rho _{1}=\beta \nu _{1}+\gamma \mu _{1},\\
\varphi (F_{1}e_{0}) &=&\varphi (F_{1})\varphi (e_{0})\Rightarrow \sigma
_{0}^{^{\prime }}=\tau _{0}^{\prime }=0,\omega _{1}=\beta \tau _{1}+\gamma
\sigma _{1}, \\
\varphi (F_{0}E_{1}) &=&\varphi (F_{0})\varphi (E_{1})\Rightarrow \varphi
(C_{0})=\left( \sigma _{0}\nu _{1}+\tau _{0}\mu _{1}\right) C_{0}, \\
\varphi (E_{0}F_{1}) &=&\varphi (E_{0})\varphi (F_{1})\Rightarrow \varphi
(C_{0})=\left( \mu _{0}\tau _{1}+\nu _{0}\sigma _{1}\right) C_{0}, \\
\varphi (E_{1}F_{0}) &=&\varphi (E_{1})\varphi (F_{0})\Rightarrow \varphi
(C_{1})=\left( \mu _{1}\tau _{0}+\nu _{1}\sigma _{0}\right) C_{1}, \\
\varphi (F_{1}E_{0}) &=&\varphi (F_{1})\varphi (E_{0})\Rightarrow \varphi
(C_{1})=\left( \mu _{0}\tau _{1}+\nu _{0}\sigma _{1}\right) C_{1}, \\
\varphi (E_{0}E_{1}) &=&\varphi (E_{0})\varphi (E_{1})\Rightarrow \mu
_{0}\nu _{1}+\nu _{0}\mu _{1}=0, \\
\varphi (F_{0}F_{1}) &=&\varphi (F_{0})\varphi (F_{1})\Rightarrow \sigma
_{0}\tau _{1}+\tau _{0}\sigma _{1}=0.
\end{eqnarray*}

2) Let $\varphi \in Aut_{I}(H_{1}^{i})$, $\psi =\varphi _{1}\kappa \in
Aut_{II}(H_{1}^{i})$, $\varphi _{1}\in Aut_{I}(H_{1}^{i})$, since $\psi
^{-1}\varphi \psi =\kappa ^{-1}\varphi _{1}^{-1}\varphi \varphi _{1}\kappa $%
, it suffices to prove that $ad\kappa $ leaves $Aut_{I}(H_{1}^{i})$ stable,
now with $\varphi $ as in \ref{34},\ref{35}, one has $\kappa ^{-1}\varphi
\kappa =\varphi ^{\prime }$, $\varphi ^{\prime }$ of the form $\ref{34}$ and $\ref{35}$, with $\beta ^{\prime }=-\delta ,\gamma ^{\prime }=-\eta ,\delta
^{\prime }=-\beta ,\eta ^{\prime }=-\gamma ,$ $\mu _{0}^{\prime }=\mu
_{1},\nu _{0}^{\prime }=\nu _{1},$ $\sigma _{0}^{\prime }=\sigma _{1},\tau
_{0}^{\prime }=\tau _{1},\mu _{1}^{\prime }=\mu _{0},\nu _{1}^{\prime }=\nu
_{0},\sigma _{1}^{\prime }=\sigma _{0},\tau _{1}^{\prime }=\tau _{0}$.

3) Follows from

\begin{eqnarray*}
M^{t}GM &=&\left( 
\begin{array}{llll}
\mu _{0} & \sigma _{0} & 0 & 0 \\ 
\nu _{0} & \tau _{0} & 0 & 0 \\ 
0 & 0 & \mu _{1} & \sigma _{1} \\ 
0 & 0 & \nu _{1} & \tau _{1}
\end{array}
\right) \left( 
\begin{array}{llll}
0 & 0 & a & b \\ 
0 & 0 & -b & d \\ 
a & -b & 0 & 0 \\ 
b & d & 0 & 0
\end{array}
\right) \left( 
\begin{array}{llll}
\mu _{0} & \nu _{0} & 0 & 0 \\ 
\sigma _{0} & \tau _{0} & 0 & 0 \\ 
0 & 0 & \mu _{1} & \nu _{1} \\ 
0 & 0 & \sigma _{1} & \tau _{1}
\end{array}
\right) \\
&=&\left( 
\begin{array}{llll}
0 & 0 & a\mu _{0}-b\sigma _{0} & b\mu _{0}+d\sigma _{0} \\ 
0 & 0 & a\nu _{0}-b\tau _{0} & -b\nu _{0}+d\tau _{0} \\ 
a\mu _{1}+b\sigma _{1} & -b\mu _{1}+d\sigma _{1} & 0 & 0 \\ 
a\nu _{1}+b\tau _{1} & -b\nu _{1}+d\tau _{1} & 0 & 0
\end{array}
\right) \left( 
\begin{array}{llll}
\mu _{0} & \nu _{0} & 0 & 0 \\ 
\sigma _{0} & \tau _{0} & 0 & 0 \\ 
0 & 0 & \mu _{1} & \nu _{1} \\ 
0 & 0 & \sigma _{1} & \tau _{1}
\end{array}
\right).
\end{eqnarray*}

Thus, the symmetric matrix $M^{t}GM$ will be of the form 
$$
\left( 
\begin{array}{llll}
0 & 0 & a^{\prime } & b^{\prime } \\ 
0 & 0 & -b^{\prime } & d \\ 
a^{\prime } & -b^{\prime } & 0 & 0 \\ 
b^{\prime } & d^{\prime } & 0 & 0
\end{array}
\right) 
$$
for $a^{\prime },b^{\prime },c^{\prime },d^{\prime }\in \ccc$, whenever the sum $X+Y$ for its entries $X,Y$ located as 
$$
\left( 
\begin{array}{llll}
. & . & . & X \\ 
. & . & Y & . \\ 
. & . & . & . \\ 
. & . & . & .
\end{array}
\right) 
$$ 
vanishes. The latter are
$$
\left\{ 
\begin{array}{c}
X=\left( a\mu _{0}-b\sigma _{0}\right) \mu _{1}+\left( b\mu _{0}+d\sigma
_{0}\right) \nu _{1}, \\ 
Y=\left( a\nu _{0}-b\tau _{0}\right) \mu _{1}+\left( b\nu _{0}+d\tau
_{0}\right) \nu _{1},
\end{array}
\right. 
$$
whenever the vanishing of
$$
X+Y=a\left( \nu _{1}\mu _{0}+\mu _{1}\nu _{0}\right) +b\left( \tau _{1}\mu
_{0}+\sigma _{1}\nu _{0}-\nu _{1}\sigma _{0}-\mu _{1}\tau _{0}\right)
+d\left( \tau _{1}\sigma _{0}+\sigma _{1}\tau _{0}\right) 
$$
expresses the constraints in $\ref{36}$.

4) Multiplying both sides of \ref{36}-(3) by $\tau _{1}\nu _{1}$ yield using 
$\ref{36}$-(1) and $\ref{36}$-(2) yields
$$
-\nu _{0}\tau _{1}D_{1}=\tau _{0}\nu _{1}D_{1}, 
$$
with $D_{1}=\mu _{1}\tau _{1}-\nu _{1}\sigma _{1}$ assumed not to vanish
thus implying relation $\ref{37}$
\edemo

\subsubsection{The group $Int(H_{1}^{i})$ of inner automorphisms of $%
H_{1}^{i}$}

\begin{proposition}
1) The element $%
h=a_{0}e_{0}+X_{0}E_{0}+Y_{0}F_{0}+c_{0}C_{0}+a_{1}e_{1}+X_{1}E_{1}+Y_{1}F_{1}+c_{1}C_{1}\in H_{1}^{i}
$ indexed by $a_{0},X_{0},Y_{0},c_{0},a_{1},X_{1},Y_{1},c_{1}\in \ccc$,
is invertible iff $a_{0}a_{1}\neq 0$. Its inverse is then 
\begin{eqnarray}
h^{-1} &=&\frac{1}{a_{0}}e_{0}-a_{0}a_{1}\left( X_{0}E_{0}+Y_{0}F_{0}\right)
+\frac{1}{a_{0}}\left( \frac{P}{a_{0}a_{1}}-\frac{1}{a_{0}}c_{0}\right) C_{0}
\label{38}, \\
&&+\frac{1}{a_{1}}e_{1}-a_{0}a_{1}\left( X_{1}E_{1}+Y_{1}F_{1}\right) +\frac{%
1}{a_{1}}\left( \frac{P}{a_{0}a_{1}}-\frac{1}{a_{1}}c_{1}\right) C_{1}, 
\nonumber
\end{eqnarray}

where $P=X_{0}Y_{1}+X_{1}Y_{0}$.

2) The matrix of $ad(h)$ then reads 
\begin{equation}
\left(
\begin{array}{cccccccc}
1 & 0 & 0 & 0 & 0 & 0 & 0 & 0 \\ 
-\frac{X_{0}}{a_{1}} & \frac{a_{0}}{a_{1}} & 0 & 0 & \frac{X_{0}}{a_{1}}&
 0 & 0&0\\ 
-\frac{Y_{0}}{a_{1}} & 0 & \frac{a_{0}}{a_{1}} & 0 & \frac{Y_{0}}{a_{1}}
 & 0 & 0 & 0 \\ 
\frac{P}{a_{0}a_{1}} & -\frac{Y_{1}}{a_{1}} & -\frac{X_{1}}{a_{1}} & 
1 & -\frac{P}{a_{0}a_{1}} & \frac{Y_{0}}{a_{0}} & \frac{X_{0}}{a_{0}} & 0 \\ 
0 & 0 & 0 & 0 & 1 & 0 & 0 & 0 \\ 
\frac{X_{1}}{a_{0}} & 0 & 0 & 0 & -\frac{X_{1}}{a_{0}} & \frac{a_{1}}{a_{0}} & 0 & 0 \\ 
\frac{Y_{1}}{a_{0}} & 0 & 0 & 0 & -\frac{Y_{1}}{a_{0}} & 0 & \frac{%
a_{1}}{a_{0}} & 0 \\ 
-\frac{P}{a_{0}a_{1}} & \frac{Y_{1}}{a_{1}} & \frac{X_{1}}{a_{1}} & 0
& \frac{P}{a_{0}a_{1}} & -\frac{Y_{0}}{a_{0}} & -\frac{X_{0}}{a_{0}} & 1
\end{array}
\right)
\label{39}
\end{equation}
in the basis $(e_{0},E_{0},F_{0},C_{0},e_{1},E_{1},F_{1},C_{1})$.

This displays the four-parametric group $Int(H_{1}^{i})\subset
Aut_{I}(H_{1}^{i})$ of inner automophisms of $H_{1}^{i}$. The relationship
with the parametrization of the full $Aut_{I}(H_{1}^{i})$ is as follows. 
\begin{equation}
\left\{ 
\begin{array}{c}
\beta =-\frac{X_{0}}{a_{1}}, \\ 
\gamma =-\frac{Y_{0}}{a_{1}}, \\ 
\delta =\frac{X_{1}}{a_{0}}, \\ 
\eta =\frac{Y_{1}}{a_{0}},
\end{array}
\right. \left\{ 
\begin{array}{c}
\mu _{0}=\frac{a_{0}}{a_{1}}, \\ 
\nu _{0}=0, \\ 
\sigma _{0}=0, \\ 
\tau _{0}=\frac{a_{0}}{a_{1}},
\end{array}
\right. \left\{ 
\begin{array}{c}
\mu _{1}=\frac{a_{1}}{a_{0}}, \\ 
\nu _{1}=0, \\ 
\sigma _{1}=0,\\ 
\tau _{1}=\frac{a_{1}}{a_{0}}.
\end{array}
\right.   \label{40}
\end{equation}
\end{proposition}

\demo
Let $h^{\prime }=a_{0}^{\prime }e_{0}+X_{0}^{\prime }E_{0}+Y_{0}^{\prime
}F_{0}+c_{0}^{\prime }C_{0}+a_{1}^{\prime }e_{1}+X_{1}^{\prime
}E_{1}+Y_{1}^{\prime }F_{1}+c_{1}^{\prime }C_{1}\in H_{1}^{i}$, we have $%
hh^{\prime }=\mathbf{1}$ iff

$$
\left\{ 
\begin{array}{c}
a_{0}a_{0}^{\prime }=1, \\ 
a_{0}X_{0}^{\prime }+a_{1}^{\prime }X_{0}=0, \\ 
a_{0}Y_{0}^{\prime }+a_{1}^{\prime }Y_{0}=0, \\ 
a_{0}c_{0}^{\prime }+c_{0}a_{0}^{\prime }+X_{0}Y_{1}^{\prime
}+Y_{0}X_{1}^{\prime }=0,
\end{array}
\right. \quad \left\{ 
\begin{array}{c}
a_{1}a_{1}^{\prime }=1, \\ 
a_{0}^{\prime }X_{1}+a_{1}X_{1}^{\prime }=0, \\ 
a_{1}Y_{1}^{\prime }+a_{0}^{\prime }Y_{1}=0, \\ 
a_{1}c_{1}^{\prime }+c_{1}a_{1}^{\prime }+X_{1}Y_{0}^{\prime
}+Y_{1}X_{0}^{\prime }=0,
\end{array}
\right. 
$$

hence $\ref{38}$ is proved. With $\lambda $ (resp. $\rho $) the regular representation (resp. antirepresentation), we have

$$
matrix\;of\;\lambda (h)=\left( 
\begin{array}{lllllllll}
& e_{0} & E_{0} & F_{0} & C_{0} & e_{1} & E_{1} & F_{1} & C_{1} \\ 
e_{0} & a_{0} & 0 & 0 & 0 & 0 & 0 & 0 & 0 \\ 
E_{0} & 0 & a_{0} & 0 & 0 & X_{0} & 0 & 0 & 0 \\ 
F_{0} & 0 & 0 & a_{0} & 0 & Y_{0} & 0 & 0 & 0 \\ 
C_{0} & c_{0} & 0 & 0 & a_{0} & 0 & Y_{0} & X_{0} & 0 \\ 
e_{1} & 0 & 0 & 0 & 0 & a_{1} & 0 & 0 & 0 \\ 
E_{1} & X_{1} & 0 & 0 & 0 & 0 & a_{1} & 0 & 0 \\ 
F_{1} & Y_{1} & 0 & 0 & 0 & 0 & 0 & a_{1} & 0 \\ 
C_{1} & 0 & Y_{1} & X_{1} & 0 & c_{1} & 0 & 0 & a_{1}
\end{array}
\right), 
$$

and

$$
matrix\;of\;\rho (h^{\prime })=\left( 
\begin{array}{lllllllll}
& e_{0} & E_{0} & F_{0} & C_{0} & e_{1} & E_{1} & F_{1} & C_{1} \\ 
e_{0} & a_{0}^{\prime } & 0 & 0 & 0 & 0 & 0 & 0 & 0 \\ 
E_{0} & X_{0}^{\prime } & a_{1}^{\prime } & 0 & 0 & 0 & 0 & 0 & 0 \\ 
F_{0} & Y_{0}^{\prime } & 0 & a_{1}^{\prime } & 0 & 0 & 0 & 0 & 0 \\ 
C_{0} & c_{0}^{\prime } & Y_{1}^{\prime } & X_{1}^{\prime } & a_{0}^{\prime }
& 0 & 0 & 0 & 0 \\ 
e_{1} & 0 & 0 & 0 & 0 & a_{1}^{\prime } & 0 & 0 & 0 \\ 
E_{1} & 0 & 0 & 0 & 0 & X_{1}^{\prime } & a_{1}^{\prime } & 0 & 0 \\ 
F_{1} & 0 & 0 & 0 & 0 & Y_{1}^{\prime } & 0 & a_{1}^{\prime } & 0 \\ 
C_{1} & 0 & 0 & 0 & 0 & c_{1}^{\prime } & Y_{0}^{\prime } & X_{0}^{\prime }
& a_{1}^{\prime }
\end{array}
\right). 
$$

The matrix product $\lambda (h)\rho (h^{\prime })$ then yields relation $\ref{39}$
\edemo

\subsubsection{The subgroup $Aut^{S}(H_{1}^{i})$ of automorphisms of $%
H_{1}^{i}$ commuting with the antipode $S$}

\begin{proposition}
The $\varphi \in Aut^{S}(H_{1}^{i})$ are of the following two types,
corresponding to type $Aut_{I}(H_{1}^{i})$ and type $Aut_{II}(H_{1}^{i})$
with $\mu ,\nu ,\sigma ,\tau \in \ccc$ such that $\mu \nu -\sigma \tau
\neq 0$.

1) The $\varphi \in Aut_{I}^{S}(H_{1}^{i})$ are as follows 
$$
\left\{ 
\begin{array}{c}
\varphi (e_{0})=e_{0}, \\ 
\varphi (E_{0})=\mu E_{0}+\nu F_{0}, \\ 
\varphi (F_{0})=\sigma E_{0}+\tau F_{0}, \\ 
\varphi (C_{0})=\left( \mu \tau -\sigma \nu \right) C_{0},
\end{array}
\right. \left\{ 
\begin{array}{c}
\varphi (e_{1})=e_{1}, \\ 
\varphi (E_{1})=\mu E_{1}-\nu F_{1}, \\ 
\varphi (F_{1})=-\sigma E_{1}+\tau F_{1}, \\ 
\varphi (C_{1})=\left( \mu \tau -\sigma \nu \right) C_{1}.
\end{array}
\right. 
$$
Their action on the generators is given by 
\begin{equation}
\left\{ 
\begin{array}{c}
\varphi (K)=K, \\ 
\varphi (E)=\mu E+\nu KF, \\ 
\varphi (F)=\sigma KE+\tau F.
\end{array}
\right.   \label{41}
\end{equation}

2) The $\varphi \in Aut_{II}^{S}(H_{1}^{i})$ are as follows 
$$
\left\{ 
\begin{array}{c}
\varphi (e_{0})=e_{1}, \\ 
\varphi (E_{0})=\mu E_{1}+\nu F_{1}, \\ 
\varphi (F_{0})=\sigma E_{1}+\tau F_{1}, \\ 
\varphi (C_{0})=\left( \mu \tau -\sigma \nu \right) C_{1},
\end{array}
\right. \left\{ 
\begin{array}{c}
\varphi (e_{1})=e_{0}, \\ 
\varphi (E_{1})=-\mu E_{0}+\nu F_{0}, \\ 
\varphi (F_{1})=\sigma E_{0}-\tau F_{0}, \\ 
\varphi (C_{1})=\left( \mu \tau -\sigma \nu \right) C_{0}.
\end{array}
\right. 
$$

Their action on the generators is given by 
\begin{equation}
\left\{ 
\begin{array}{c}
\varphi (K)=-K, \\ 
\varphi (E)=-\mu KE+\nu F, \\ 
\varphi (F)=\sigma E-\tau KF.
\end{array}
\right.   \label{42}
\end{equation}
\end{proposition}

\demo
We confer $\ref{13}$ with the action of $S$ which we recall.

\begin{center}
\begin{tabular}{|c|c|c|c|c|c|c|c|c|}
\hline
$a$ & $e_{0}$ & $E_{0}$ & $F_{0}$ & $C_{0}$ & $e_{1}$ & $E_{1}$ & $F_{1}$ & $C_{1}$ \\
\hline 
$S(a)$ & $e_{0}$ & $E_{1} $&$ -F_{1} $& $C_{0} $&$ e_{1}$ & $-E_{0}$ & $F_{0}$ & $C_{1}$\\
\hline
\end{tabular}
\end{center}

1) The requirement that the action of $\varphi \circ S$ and $S\circ \varphi $
be the same,

-on $e_{0}$ yields $\beta =\gamma =\delta =\eta =0$,

-on $F_{0}$ or $F_{1}$ yields $\sigma _{0}=-\sigma _{1}$ and $\tau _{0}=\tau
_{1}$,

-on $C_{0}$ is automatic.

The expression $\ref{40}$ immediately follows from $\ref{39}$ and from the
facts that $E_{0}+E_{1}=E,F_{0}+F_{1}=F,E_{0}-E_{1}=KE,F_{0}-F_{1}=KF$.

2) The requirement that the action of $\varphi \circ S$ and $S\circ \varphi $
be the same,

-on $e_{1}$ yields $\beta =\gamma =\delta =\eta =0$,

-on $F_{1}$ or $F_{0}$ yields $\sigma _{0}=-\sigma _{1}$ and $\tau _{0}=\tau
_{1}$,

-on $C_{1}$ is automatic,
\edemo

\subsubsection{The Hopf automorphisms of $H_{1}^{i}$}

\begin{corollary}
Each element of $Aut_{I}^{S}(H_{1}^{i})$ and none of $Aut_{II}^{S}(H_{1}^{i})$ is a Hopf automorphism.
\end{corollary}

\demo
We check that $\varphi \in Aut_{I}^{S}(H_{1}^{i})$ is coalgebra morphism. By
the multiplicativity of $\Delta$, it suffices to address the generators. We have, by $\ref{41}$,

\begin{eqnarray*}
\Delta (\varphi (K)) &=&\Delta (K)=K\otimes K=\varphi (K)\otimes \varphi (K),
\\
\Delta (\varphi (E)) &=&\Delta \left( \mu E+\nu KF\right) =\mu \left(
E\otimes \mathbf{1}+K\otimes E\right) +\nu \left( K\otimes K\right) \left(
F\otimes K^{-1}+\mathbf{1}\otimes F\right), \\
&=&\left( \mu E+\nu KF\right) \otimes \mathbf{1+}K\otimes \left( \mu E+\nu
KF\right) =\varphi (E)\otimes \varphi (\mathbf{1})+\varphi (K)\otimes
\varphi \left( F\right), \\
\Delta (\varphi (F)) &=&\Delta \left( \sigma KE+\tau F\right) =\sigma \left(
K\otimes K\right) \left( E\otimes \mathbf{1}+K\otimes E\right) +\tau \left(
F\otimes K^{-1}+\mathbf{1}\otimes F\right), \\
&=&\left( \sigma KE+\tau F\right) \otimes K^{-1}+\mathbf{1}\otimes \left(
\sigma KE+\tau F\right) =\varphi (F)\otimes \varphi (K^{-1})+\varphi (%
\mathbf{1})\otimes \varphi (F),
\end{eqnarray*}
whereas, for $\varphi \in Aut_{II}^{S}(H_{1}^{i})$, we have by \ref{42}, 
$$
\Delta (\varphi (K))=-\Delta (K)=-K\otimes K\neq \varphi (K)\otimes \varphi
(K)=K\otimes K.
$$
\edemo

\begin{lemma}
The identity on the generators $K,E,F$ extends uniquely to a Hopf $*$%
-operation $I$ of $H_{1}^{i}$, whose action is given as follows. 
\begin{center}
\begin{tabular}{|c|c|c|c|c|c|c|c|c|}
\hline
$a$ & $e_{0}$ & $E_{0}$ & $F_{0} $& $C_{0} $& $e_{1}$ & $E_{1}$ & $F_{1}$ & $C_{1}$ \\ 
\hline
$I(a)$ & $e_{0}$ & $E_{1}$ & $F_{1}$ & $C_{0}$ &$ e_{1} $& $E_{0}$ & $F_{0}$ & $C_{1}$\\
\hline
\end{tabular}
\end{center}
\end{lemma}

\demo
The defining $\ref{11}$ and $\ref{12}$ are obviously respected, as well as the
definition relations of the Hopf structure stated on the generators.
\edemo

\begin{proposition}
1) The semi-Hopf $*$-operations of $H_{1}^{i}$ are the of the following two types,
corresponding to the above type $I$ and type $II$ automorphisms.

Type $I$: $\Gamma =I\circ \varphi ,\varphi \in Aut_{I}^{S}(H_{1}^{i})$
is given by 
\begin{equation}
\left\{ 
\begin{array}{c}
\Gamma (e_{0})=e_{0}, \\ 
\Gamma (E_{0})=\alpha E_{1}+\beta F_{1}, \\ 
\Gamma (F_{0})=\gamma E_{1}+\delta F_{1}, \\ 
\Gamma (C_{0})=\lambda C_{0},
\end{array}
\right. \left\{ 
\begin{array}{c}
\Gamma (e_{1})=e_{1}, \\ 
\Gamma (E_{1})=\alpha E_{0}-\beta F_{0}, \\ 
\Gamma (F_{1})=-\gamma E_{0}+\delta F_{0}, \\ 
\Gamma (C_{1})=\lambda C_{1},
\end{array}
\right.   \label{44}
\end{equation}

where 
\begin{equation}
\left\{ 
\begin{array}{l}
\alpha =ae^{i\phi }, \\ 
\beta =\pm be^{i\frac{(\phi +\psi )}{2}}, \\ 
\gamma =\pm ce^{i\frac{(\phi +\psi )}{2}}, \\ 
\delta =ae^{i\psi }, \\ 
\lambda =\left( \alpha \delta -\beta \gamma \right) =e^{i(\phi +\psi )},
\end{array}
\right. \label{45}
\end{equation}
with $a,b,c\geq 0$ fulfilling $a^{2}-bc=1$.\\ 
Type $II$: $\Gamma =I\circ \varphi ,\varphi \in Aut_{II}^{S}(H_{1}^{i})$ is given by
\begin{equation}
\left\{ 
\begin{array}{c}
\Gamma (e_{0})=e_{1}, \\ 
\Gamma (E_{0})=\alpha E_{0}+\beta F_{0}, \\ 
\Gamma (F_{0})=\gamma E_{0}+\delta F_{0}, \\ 
\Gamma (C_{0})=\lambda C_{1},
\end{array}
\right. \left\{ 
\begin{array}{c}
\Gamma (e_{1})=e_{0}, \\ 
\Gamma (E_{1})=-\alpha E_{1}+\beta F_{1}, \\ 
\Gamma (F_{1})=\gamma E_{1}-\delta F_{1}, \\ 
\Gamma (C_{1})=\lambda C_{0},
\end{array}
\right.   \label{46}
\end{equation}
where 
\begin{equation}
\left\{ 
\begin{array}{l}
\alpha =ae^{i\phi }, \\ 
\beta =\pm ibe^{i\frac{(\phi +\psi )}{2}}, \\ 
\gamma =\pm ice^{i\frac{(\phi +\psi )}{2}}, \\ 
\delta =ae^{i\psi }, \\ 
\lambda =-\left( \alpha \delta -\beta \gamma \right) =-e^{i(\phi +\psi )},
\end{array}
\right.\label{47}
\end{equation}
with $a>0,b,c\geq 0$ fulfilling $a^{2}+bc=1$.\\
2) The $*$-operations characterized under ($I$) above are in fact the Hopf $*
$-operations of $H_{1}^{i}$. Indeed they all fulfill $\Delta (\Gamma
(a))=\Gamma \otimes \Gamma (\Delta (a)),a\in H_{1}^{i}$, whilst this is the
case for none of the $*$-operations ($II$). One has thus a four-parameter
family of Hopf $*$-operations belonging to the same orbit of right action of
Hopf homomorphisms. Observe that the $*$-operation $I$ is obtained by
making in ($I$) the choice $\alpha =\delta =1,\beta =\gamma =0$ whilst the
choice $\alpha =\delta =0,\beta =-i,\gamma =i$ yields the $*$-operation 
$$
\left\{ 
\begin{array}{c}
e_{0}\rightarrow e_{0}, \\ 
E_{0}\rightarrow -iF_{1}, \\ 
F_{0}\rightarrow iE_{1}, \\ 
C_{0}\rightarrow -C_{0},
\end{array}
\right. \left\{ 
\begin{array}{c}
e_{1}\rightarrow e_{1}, \\ 
E_{1}\rightarrow iF_{0}, \\ 
F_{1}\rightarrow -iE_{0}, \\ 
C_{1}\rightarrow -C_{1},
\end{array}
\right. \left\{ 
\begin{array}{c}
K\rightarrow K, \\ 
E\rightarrow iKF, \\ 
F\rightarrow iEK^{-1}.
\end{array}
\right. 
$$
\end{proposition}

\demo

1) We seek the semi-Hopf $*$-operations as the composition products $I\circ
\varphi $, $\varphi \in Aut^{S}(H_{1}^{i})$, which are involutions. With $*$
indicating complex conjugation, we have 

-for $\varphi \in Aut_{I}^{S}(H_{1}^{i})$ iteration of $\ref{44}$ will yield
the identity operation iff $M\bar{M}=\mathbf{1}$, $M$ the matrix $\left( 
\begin{array}{ll}
\mu & \nu \\ 
\sigma & \tau
\end{array}
\right) $, and $\lambda \bar{\lambda}=1$. Setting $\mu ^{*}=\alpha ,\nu
^{*}=\beta ,\sigma ^{*}=\gamma ,\tau ^{*}=\delta $, the first condition
yields

$$
\left\{ 
\begin{array}{c}
1)\alpha \alpha ^{*}-\beta \gamma ^{*}=1, \\ 
2)\delta \delta ^{*}-\beta \gamma ^{*}=1,
\end{array}
\right. \left\{ 
\begin{array}{c}
3)\alpha \beta ^{*}=\beta \delta ^{*}, \\ 
4)\gamma \alpha ^{*}=\delta \gamma ^{*},
\end{array}
\right. 
$$

expressed by \ref{45} (observe that $a\neq 0$ and that in the case $b=c=0$
the phases of $\beta $ and $\gamma $ are arbitrary); the second condition is
then automatic,

-for $\varphi \in Aut_{II}^{S}(H_{1}^{i})$ iteration of $\ref{46}$ will yield
the identity operation iff $M\bar{M}=\mathbf{1}$, $M$ the matrix $\left( 
\begin{array}{ll}
\mu & -\nu \\ 
-\sigma & \tau
\end{array}
\right) $, and $\lambda \bar{\lambda}=1$. Setting $\mu ^{*}=\alpha ,\nu
^{*}=\beta ,\sigma ^{*}=\gamma ,\tau ^{*}=\delta $, the first condition
yields

$$
\left\{ 
\begin{array}{c}
1)\alpha \alpha ^{*}+\beta \gamma ^{*}=1, \\ 
2)\delta \delta ^{*}+\beta \gamma ^{*}=1,
\end{array}
\right. \left\{ 
\begin{array}{c}
3)\alpha \beta ^{*}+\beta \delta ^{*}=0, \\ 
4)\gamma \alpha ^{*}+\delta \gamma ^{*}=0,
\end{array}
\right. 
$$

expressed by $\ref{47}$ (observe that $a\neq 0$ and that in the case $b=c=0$
the phases of $\beta $ and $\gamma $ are arbitrary); the second condition is
then automatic.

2)Follows from 5.2.2. or can be checked analogously.
\edemo


\begin{thebibliography}{CIKS}
\bibitem{AC1}  A. Connes {\it Brisure spontan\'ee de sym\'etrie et g\'eom\'etrie du point de vue spectral}, S\'eminaire Bourbaki, 48\`emem ann\'ee, {\bf 816} (1996)

\bibitem{CC1}  A. Chamseddine and A. Connes {\it The spectral action principle}, hep-th/9606001, to appear in Comm. in Math, Phys. (1996)

\bibitem{CIKS}  L. Carminati, B. Iochum, D.Kastler and T. Sch\"ucker {\it On connes' new principle of general relativity: can spinors hear the forces of space-time?}, hep-th/9612228 (1996) 

\bibitem{DK1}  D. Kastler {\it Regular and adjoint representation of $SL_{q}(2)$ at third root of unit}, CPT internal report (1995)

\bibitem{DK2}  D. Kastler {\it Introduction \`a l'\'electrodynamique quantique}, Dunod, Paris (1960)

\bibitem{CASS}  C. Kassel {\it Quantum Groups}, Springer, Berlin (1995)

\end{thebibliography}
\end{document}